\documentclass[a4paper,11pt,twoside]{article}
\usepackage{amssymb,amsmath,latexsym,geometry,lineno,fancyhdr,hyperref,titletoc}
\usepackage{amsthm}

\usepackage{color}

\usepackage{import}
\usepackage{xifthen}
\usepackage{pdfpages}
\usepackage{amsmath}
\usepackage{mathrsfs}
\contentsmargin{0pt}

\dottedcontents{section}[2em]{\vspace{-1mm}\small}{2em}{0pt}
\dottedcontents{subsection}[5em]{\vspace{-1mm}\small}{3em}{5pt}

\hypersetup{colorlinks,linkcolor= blue,citecolor=blue}

\newtheorem{theorem}{Theorem}[section]

\newtheorem{lemma}{Lemma}[section]

\newtheorem{proposition}{Proposition}[section]

\newtheorem{remark}{Remark}[section]

\newtheorem{claim}{Claim}[section]

\newtheorem{corollary}{Corollary}[section]

\newtheorem{definition}{Definition}[section]

\numberwithin{equation}{section}
\allowdisplaybreaks[0]

\def\Vs{\vskip8pt}\def\vs{\vskip4pt}

\begin{document}
\begin{sloppypar}
\title{{\bf\Large Global bifurcations of nodal solutions for coupled elliptic equations}}
\author{\\
{\textbf{\normalsize Haoyu Li, Ol\'impio Hiroshi Miyagaki}}\\
{\it\small Departamento de Matem\'atica, Universidade Federal de S\~{a}o Carlos}\\
{\it\small S\~{a}o Carlos-SP, 13565-905, Brazil}\\
{\textbf{\normalsize Zhi-Qiang Wang}}\\
{\it\small Department of Mathematics and Statistics, Utah State University}\\
{\it\small Logan, UT 84322, USA}\\
}
\date{}
\maketitle
{\bf\normalsize Abstract.} {\small
We investigate the global bifurcation structure of the radial nodal solutions to the coupled elliptic equations
\begin{equation}
    \left\{
   \begin{array}{lr}
     -{\Delta}u+u=u^3+\beta uv^2\mbox{ in }B_1 ,\nonumber\\
     -{\Delta}v+v=v^3+\beta u^2v\mbox{ in }B_1 ,\nonumber\\
     u,v\in H_{0,r}^1(B_1).\nonumber
   \end{array}
   \right.
\end{equation}
Here $B_1$ is a unit ball in $\mathbb{R}^3$ and $\beta\in\mathbb{R}$ the coupling constant is used as bifurcation parameter.
For each $k$, the unique pair of nodal solutions $\pm w_k$ with exactly $k-1$ zeroes to the scalar field equation $-\Delta w + w=w^3$ generate exactly four synchronized solution curves and exactly four semi-trivial solution curves to the above system. We obtain a fairly complete global bifurcation structure of all
bifurcating branches emanating from these eight solution curves of the system, and show that for different $k$ these bifurcation structures are disjoint.
We obtain exact and distinct nodal information for each of the bifurcating branches, thus providing a fairly complete characterization of nodal solutions of the system in terms of the coupling.}

\footnotetext{e-mail of HL: hyli1994@hotmail.com, e-mail of OHM: olimpio@ufscar.br, e-mail of ZQW: zhi-qiang.wang@usu.edu}



\pagestyle{fancy}
\fancyhead{} 
\fancyfoot{} 
\renewcommand{\headrulewidth}{0pt}
\renewcommand{\footrulewidth}{0pt}
\fancyhead[CE]{ \textsc{Haoyu Li \& Ol\'impio Hiroshi Miyagaki \& Zhi-Qiang Wang}}
\fancyhead[CO]{ \textsc{Global bifurcation of nodal solution for coupled elliptic equation}}
\fancyfoot[C]{\thepage}

\tableofcontents

\section{Introduction}
\subsection{Introduction}
We investigate the global bifurcation structure for the following coupled nonlinear elliptic equations:
\begin{equation}\label{e:001}
    \left\{
   \begin{array}{lr}
     -{\Delta}u+u=u^3+\beta uv^2\mbox{ in }B_1 ,\\
     -{\Delta}v+v=v^3+\beta u^2v\mbox{ in }B_1 ,\\
     u,v\in H_{0,r}^1(B_1).
   \end{array}
   \right.
\end{equation}
Here, the domain $B_1$ is the three dimensional unit ball centered at the origin and $\beta\in\mathbb{R}$ is considered the bifurcation parameter for the system.
Such kind of systems and related systems arise in many physical phenomenon, especially in nonlinear optics \cite{AkhmedievAnkiewicz1999} for the Kerr-like photorefractive media. Moreover, one can find Problem (\ref{e:001}) in other physical problems, such as in the Hartree-Fock theory of a double condensate \cite{ErsyGreeneBurke1997}. The constant $\beta\in\mathbb{R}$ is called the coupling constant which describes the interaction between interspecies. It is called attractive if $\beta$ is positive and repulsive if $\beta$ is negative.

Problem (\ref{e:001}) and related ones have attracted extensive interests in recent years due to their rich and complex solution structure both in physics background and mathematical analysis, following the initial work of Lin-Wei \cite{LinWei2005} there have been extensive works cf. \cite{AmbrosettiColorado2007,LiuWang2008,MaiaMontefuscoPellacci2008,Hirano2009,TerraciniVerzini2009,ChenZou2011,PengWang2013,ClappPistoia2018} and the references therein.
Many of these works have been done from the point of views of variational methods and bifurcation analysis.

It is known that the positive solutions to the subcritical scalar field equation $-\Delta u+u=u^p$ possess an a priori estimate, radial symmetry up to translations, and uniqueness in the presence of spatial symmetry \cite{GidasSpruck1981, GidasNiNirenberg1979,Kwong1989,Tanaka2016}. An analogue boundedness result for Problem (\ref{e:001}) can be proved via Liouville-type theorem (see \cite{DancerWeiWeth2010,QuittnerBook}) if $\beta>-1$. However, this result is no longer valid for Problem (\ref{e:001}) with $\beta\leq-1$. It is proved in \cite{WeiWeth2008,DancerWeiWeth2010} that if $\beta\leq-1$ Problem (\ref{e:001}) admits an unbounded sequence of positive solutions, which is significantly different from the counterpart of the scalar field equation. This multiplicity result of positive solutions was in \cite{TianWang2011} generalized to the $N$-coupled system case.  In \cite{BartschDancerWang2010}, the authors gave a global description of positive solutions via a bifurcation analysis \cite{CrandallRabinowitz1971,Rabinowitz1971}. It is shown in \cite{BartschDancerWang2010} that there exists a sequence of bifurcation parameters when $\beta\downarrow-1$ and each of the bifurcation branches is unbounded in $\mathbb{R}\times H_r^1(\mathbb{R}^N)\times H_r^1(\mathbb{R}^N)$ and bounded on any compact interval of $\beta$-axis. From these results we can see that the coupling constant $\beta$ plays an important and subtle role in the structure of solutions when $\beta$ varies from attractive to the repulsive regimes.

The aim of this paper is to investigate the global bifurcation structure of radial sign-changing solutions to (\ref{e:001}). More precisely we aim at characterizing all radial nodal solutions in terms of the coupling constant $\beta$ by making use of the abstract methods as well as some concrete approaches as in \cite{CrandallRabinowitz1971,BartschWangWei2007,BartschDancerWang2010} for the radial nodal solutions of Problem (\ref{e:001}).
Before we state the main theorems, let us review some of the known results on Problem (\ref{e:001}).

There are several ways to study the radial sign-changing solutions, such as the invariant set of descending flow \cite{LiuLiuWang2015,LiWang2021}, gluing on Nehari set \cite{MaiaMontefuscoPellacci2008,LiuWang2019}, etc. In \cite{LiuLiuWang2015}, the authors proved that if $\beta\leq0$ there are infinitely many solutions to Problem (\ref{e:001}) with $u$ positive and $v$ sign-changing. This shows that the set of sign-changing solutions is more complicated than the set of positive solutions. Furthermore, if $\beta\leq0$ for any $P,Q\in\mathbb{N}$ it is shown in \cite{LiuWang2019,LiWang2021} there exists a solution $(u,v)$ such that $u$ changes its sign exactly $P$ times and $v$ changes its sign exactly $Q$ times. A multiplicity result of this type is also valid. To be precise, for any $P\in\mathbb{N}$, if $\beta\leq-1$, \cite{LiWang2021} shows that there is an unbounded sequence of solutions $(u_n,v_n)$ such that both of $u_n$ and $v_n$ change their signs exactly $P$ times.

The above discussions motivate us to study sign-changing solutions from a bifurcation point of view, which would give a more complete global description of all solutions in terms of the coupling constant $\beta$. This will be the main theme of this paper and our main result will be stated in Theorems \ref{t:bifurcation} and \ref{t:bifurcation2} below. We follow the approach of \cite{BartschDancerWang2010} to analyze the local and global bifurcations both synchronized branches and semi-trivial branches. The major difficulty here is that there are infinitely many synchronized branches and infinitely many semi-trivial branches, while for nonnegative solutions there is exactly one synchronized branch and exactly two semi-trivial branches (\cite{BartschDancerWang2010}). We not only need to consider the infinitely many bifurcating branches emanating from each synchronized branch or from each semi-trivial branch, we also need to investigate the relations between the bifurcating branches out of different synchronized or semi-trivial branches.
Our results can be roughly described as follows. For each integer $k\geq 1$, let $w_k$ be the unique nodal solution $w_k$ with exactly $k-1$ zeroes and $w_k(0)>0$ of $-\Delta w+w=w^3$. Then there are exactly four synchronized solution curves (for $\beta >-1$) and exactly four semi-trivial solution curves (for $\beta \in (-\infty, \infty)$) all generated by this $w_k$, and these eight curves form a global tree-trunk structure being connected from each other exactly at $\beta=1$ by a circle set of solutions. Our results show that for each such structure generated by $w_k$ there are infinitely many bifurcating branches extended to infinity in $\beta$ and infinitely many bifurcating branches extended to negative infinity in $\beta$, forming a global bifurcating tree structure.
We also show that these bifurcating tree structures are disjoint from each other for different $k$.
For technical ingredients, to study the bifurcation parameters and the behaviour of branches, we give some new Liouville-type theorems and a Morse index estimate for the scalar field equation Theorems \ref{t:Liouville}, \ref{t:Liouville2} and Proposition \ref{p:MorseIndex}.
For the Liouville-type theorems, for a general class of coupled elliptic systems, Quittner in \cite{Quittner2021} proved the nonexistence of solutions under the assumption of a finite number of nodes of the radial functions $u$, $v$, $u+v$ and $u-v$. We consider here Problem (\ref{e:001}), which is a special case of the problem considered by Quittner. And the assumption about the number of nodes of the functions $u$ and $v$ is removed. To be precise, we recover the nonexistence of the solution $(u,v)$ by restricting the number of nodes of $u+v$ and of $u-v$ only.

Our method and results have some immediate applications. First, we state Theorem \ref{t:nonexistencebeta=3} concerning the nonexistence of radial solution in terms of nodal number. This refines \cite{LiWang2021,LiuWang2019,MaiaMontefuscoPellacci2008} from the view point of nonexistence. The second application is Theorem \ref{t:comparison} which gives the exact number of nodes of the difference of two radial solutions of the scalar field equation. In the classical result, this number can only be bounded from above (\cite{IshiwataLi}).

\subsection{Main results: Bifurcation theorems and a nonexistence result}
First we need to fix some notations.
The energy functional of Problem (\ref{e:001}) is the $C^2$-functional $H_{0,r}^1(B_1)\times H_{0,r}^1(B_1)$ defined as
\begin{align}\label{EnergyFunctional}
J_\beta(u,v)=\frac{1}{2}\int_{B_1}|\nabla u|^2+|u|^2+|\nabla v|^2+|v|^2-\frac{1}{4}\int_{B_1}u^4+ v^4+2\beta u^2v^2.
\end{align}
Throughout the paper, we denote $H_{0,r}^1(B_1)$ the subspace of $H_{0}^1(B_1)$ consisting with radial functions.

Let us begin with the following introduction of synchronized branches and semi-trivial branches of solutions.
\begin{itemize}
  \item Let $w_k\in H_{0,r}^1(B_1)$ be the unique and non-degenerate solution to $-\Delta w+w=w^3$ with $w_k(0)>0$, which changes its sign exactly $k-1$ times for $k=1,2,\cdots$. The existence, uniqueness and the non-degeneracy are guaranteed by \cite{WillemBook1996,Tanaka2016,LiMiyagaki}, respectively.
  \item Denote $u_{k,\beta}=v_{k,\beta}=\frac{w_k}{\sqrt{1+\beta}}$ for $\beta\in(-1,+\infty)$.
\end{itemize}
We note that $(u_{k,\beta},v_{k,\beta})$ solves (\ref{e:001}) with the corresponding $\beta$.
One of the four synchronized solution curves generated by $w_k$ is defined as
\begin{align}
\mathcal{T}^1_{k}=\big\{(\beta,u_{k,\beta},v_{k,\beta})|\beta\in(-1,+\infty)\big\}.
\end{align}
Due to symmetry of the system, we also have three other synchronized solution curves generated by $w_k$
\begin{align}
\mathcal{T}^2_{k}=\big\{(\beta,-u_{k,\beta},v_{k,\beta})|\beta\in(-1,+\infty)\big\},\nonumber
\end{align}
\begin{align}
\mathcal{T}^3_{k}=\big\{(\beta,u_{k,\beta},-v_{k,\beta})|\beta\in(-1,+\infty)\big\},\nonumber
\end{align}
\begin{align}
\mathcal{T}^4_{k}=\big\{(\beta,-u_{k,\beta},-v_{k,\beta})|\beta\in(-1,+\infty)\big\}.\nonumber
\end{align}
We need the concept of nodal number for radial functions.
\begin{definition}\label{def:NodalNumber}
For any radial domain $\Omega$ and a continuous radial function $u:\Omega\to\mathbb{R}$, the nodal number of $u$ is defined as the maximum number $k\in\mathbb{N}$ such that there exists a sequence of positive numbers $x_0,x_1,\cdots,x_k\in(0,+\infty)$ such that
$\partial B_{x_i}(0)\subset\Omega$ for any $i=0,\cdots,k$;
  and $u(x)|_{|x|=x_{i-1}}\cdot u(x)|_{|x|=x_{i}}<0$ for any $i=1,\cdots,k$.
\end{definition}

Note that under the above definition, we have $n(w_k)=k-1$.

Define the set of all solutions
\begin{align}\label{SolutionSet}
\mathcal{S}:=\big\{(\beta,u,v)\in\mathbb{R}\times H_{0,r}^1(B_1)\times H_{0,r}^1(B_1)|(\beta,u,v)\mbox{ solves Problem (\ref{e:001})}\big\}.
\end{align}
For $l=1,2,3,4$, define the set of solutions outside of $\mathcal{T}^l_k$
\begin{align}
\mathcal{S}^l_k:=\big\{(\beta,u,v)\in\mathbb{R}\times H_{0,r}^1(B_1)\times H_{0,r}^1(B_1)\backslash \mathcal{T}^l_k|(\beta,u,v)\mbox{ solves Problem (\ref{e:001})}\big\}.\nonumber
\end{align}
Finally we define the four semi-trivial solution branches generated by $w_k$
\begin{align}
\mathcal{ST}_k^1:=\Big\{\Big(\beta,w_k,0\Big)\in\mathbb{R}\times H_{0,r}^1(B_1) \times H_{0,r}^1(B_1)\Big|\beta\in\mathbb{R}\Big\},\nonumber
\end{align}
\begin{align}
\mathcal{ST}_k^2:=\Big\{\Big(\beta,0,w_k\Big)\in\mathbb{R}\times H_{0,r}^1(B_1) \times H_{0,r}^1(B_1)\Big|\beta\in\mathbb{R}\Big\},\nonumber
\end{align}
\begin{align}
\mathcal{ST}_k^3:=\Big\{\Big(\beta,-w_k,0\Big)\in\mathbb{R}\times H_{0,r}^1(B_1) \times H_{0,r}^1(B_1)\Big|\beta\in\mathbb{R}\Big\},\nonumber
\end{align}
\begin{align}
\mathcal{ST}_k^4:=\Big\{\Big(\beta,0,-w_k\Big)\in\mathbb{R}\times H_{0,r}^1(B_1) \times H_{0,r}^1(B_1)\Big|\beta\in\mathbb{R}\Big\}.\nonumber
\end{align}
For $l=1,2,3,4$ define the set of solutions outside of $\mathcal{ST}^l_k$
\begin{align}
\mathcal{O}^l_k:=\big\{(\beta,u,v)\in\mathbb{R}\times H_{0,r}^1(B_1) & \times H_{0,r}^1(B_1)\backslash \mathcal{ST}^l_k|(\beta,u,v)\mbox{ solves Problem (\ref{e:001})}\big\}.\nonumber
\end{align}

\begin{definition}
For any $l=1,2,3,4$,
a parameter $\beta_0\in(-1,+\infty)$ is said to be a \emph{bifurcation parameter} of $\mathcal{T}^1_k$ if there exist a sequence $(\beta_{0,n},u_{0,n},v_{0,n})\in\mathcal{S}^l_k$ such that $(\beta_{0,n},u_{0,n},v_{0,n})\to(\beta_0,u_{\beta_0},v_{\beta_0})\in\mathcal{T}^l_k$ as $n\to\infty$. A bifurcation parameter $\beta_0$ is called a global bifurcation parameter if there exists a connected subset $\mathcal{U}$ of $\mathcal{S}^l_k$ with $(\beta_0,u_{\beta_0},v_{\beta_0})\in\overline{\mathcal{U}}$ such that either $\mathcal{U}\cup{(\beta_0,u_{\beta_0},v_{\beta_0})}$ is unbounded in $\mathbb{R}\times H_{0,r}^1(B_1)\times H_{0,r}^1(B_1)$ or $\overline{\mathcal{U}}\cap\mathcal{T}^l_k\backslash\{(\beta_0,u_{\beta_0},v_{\beta_0})\}\neq\emptyset$.
\end{definition}

\begin{remark}
This definition is related to Rabinowitz's global bifurcation theorem \cite{Rabinowitz1971}.
\end{remark}

The bifurcation parameters and the global bifurcation parameters for the four semi-trivial solution branches $\mathcal{ST}^1_k$, $\mathcal{ST}_k^2$, $\mathcal{ST}_k^3$ and $\mathcal{ST}_k^4$ can be defined in the same manner.

We use
$\mathbb{N}_+$ and $\mathbb N$ for the set of positive and nonnegative integers, and $\mbox{Proj}_\beta:\mathbb{R}\times H_{0,r}^1(B_1)\times H_{0,r}^1(B_1)\to \mathbb{R}$ for the projection $\mbox{Proj}_\beta(\beta,u,v)=\beta$.
Here are our main results.
\begin{theorem}\label{t:bifurcation}
Assume that $k\in\mathbb{N}_+$.
The bifurcation parameters of $\mathcal{T}^1_k$ form a decreasing sequence $\{\beta_{k,i}\}_{i=1}^\infty$ such that
\begin{itemize}
  \item [$(1).$] $-1<\cdots<\beta_{k,k+2}<\beta_{k,k+1}<0<1=\beta_{k,k} <\cdots<\beta_{k,1}<3$ with $\beta_{k,i}\to-1$ as $i\to\infty$;
  \item [$(2).$] for any $i\in\mathbb{N}_+\backslash\{k\}$, there is an unbounded connected set of $\mathcal{S}^1_k$, denoted by $\mathcal{U}^{(1)}_{k,i}$, such that $\overline{\mathcal{U}^{(1)}_{k,i}}\cap\mathcal{T}^1_k =\{(\beta_{k,i},u_{\beta_{k,i}},v_{\beta_{k,i}})\}$;
  \item [$(3).$] for any $i=k+1,\cdots$, $(-\infty,\beta_{k,i})\subset\mbox{Proj}_\beta(\mathcal{U}^{(1)}_{k,i}) \subset(-\infty,0)$. And for any $(\beta,u,v)\in\mathcal{U}^{(1)}_{k,i}$, we get $n(u)=n(v)=n(u+v)=k-1$, $n(u-v)=i-1$ and $u(0),v(0)>0$;
  \item [$(4).$] for any $i=1,\cdots,k-1$, $(\beta_{k,i},\infty)\subset\mbox{Proj}_\beta(\mathcal{U}^{(1)}_{k,i})\subset(1,+\infty)$.  And for any $(\beta,u,v)\in\mathcal{U}^{(1)}_{k,i}$, we get $n(u)=n(v)=n(u+v)=k-1$, $n(u-v)=i-1$ and $u(0),v(0)>0$;
  \item [$(5).$] the branch $\mathcal{U}^{(1)}_{k,k}\subset\mathcal{S}^1_k$ contains a circle of solutions
  \begin{align}
  \Big\{s_\theta=(1,\cos\theta  w_k,\sin\theta w_k) \Big|\theta\in[0,\frac{\pi}{4})\cup(\frac{\pi}{4},2\pi)\Big\}\nonumber
  \end{align}
  Moreover, $\overline{\mathcal{U}^{(1)}_{k,k}} \cap\mathcal{T}^1_k=\{(\beta_{k,k},u_{\beta_{k,k}},v_{\beta_{k,k}})\}$.
\end{itemize}
\end{theorem}

\begin{remark}
Analogue results hold for bifurcating branches emanating out of $\mathcal{T}_k^2$, $\mathcal{T}_k^3$ and $\mathcal{T}_k^4$. It is easy to verify that the corresponding bifurcation parameters are also given by $\beta_{k,i}$ due to symmetry of the system.
\end{remark}

\begin{center}
\begin{figure}[htbp]
\centering
\includegraphics[width=0.9\textwidth]{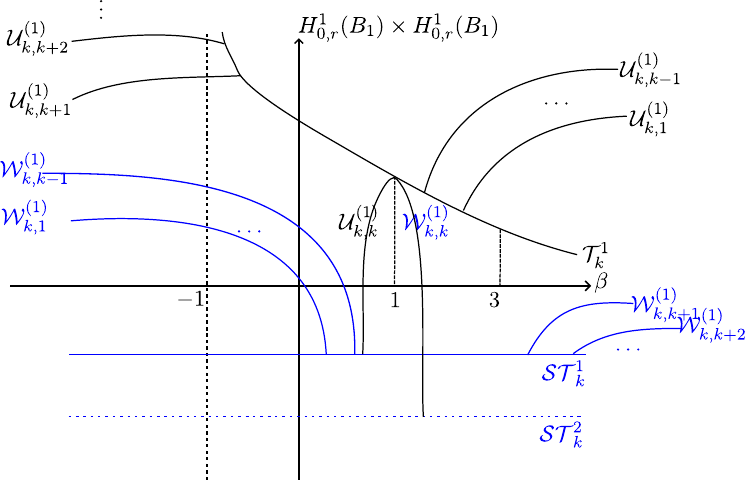}
\caption{Bifurcations emanating from $\mathcal{T}^1_k$ and $\mathcal{ST}_k^1$}
\label{fig:drawing-1-1}
\end{figure}
\end{center}

Next we consider the bifurcations on the semi-trivial branches $\mathcal{ST}_k^{l}$ for $l=1,2,3,4$, and due to symmetry we only consider $\mathcal{ST}_k^{1}$.
Let us define the numbers
\begin{align}
\widetilde{\beta}_{k,i}:=\frac{3-\beta_{k,i}}{1+\beta_{k,i}}\nonumber
\end{align}
for $k,i\in\mathbb{N}_+$.
Here, $\beta_{k,i}$ are the bifurcation parameters of $\mathcal{T}^1_k$.
\begin{theorem}\label{t:bifurcation2}
The numbers $(\widetilde{\beta}_{k,i})_{i=1}^\infty$ are the bifurcation parameters of $\mathcal{ST}^1_k$ which forms an increasing sequence such that
\begin{itemize}
  \item [$(1).$] $0<\widetilde{\beta}_{k,1}<\cdots<\widetilde{\beta}_{k,k-1}<\widetilde{\beta}_{k,k} =1<3<\widetilde{\beta}_{k,k+1}<\widetilde{\beta}_{k,k+2}<\cdots$ with $\widetilde{\beta}_{k,i}\to+\infty$ as $i\to\infty$;
  \item [$(2).$] for any $i\in\mathbb{N}_+\backslash\{k\}$, there is an unbounded connected set of $\mathcal{O}^1_k$, denoted by $\mathcal{W}^{(1)}_{k,i}$, such that $\overline{\mathcal{W}^{(1)}_{k,i}}\cap\mathcal{ST}^1_k=\{(\widetilde{\beta}^{(1)}_{k,i},w_k,0)\}$;
  \item [$(3).$] for any $i=k+1,\cdots$, $(\widetilde{\beta}_{k,i},+\infty)\subset\mbox{Proj}_\beta(\mathcal{W}^{(1)}_{k,i})\subset(3,+\infty)$. And for any $(\beta,u,v)\in\mathcal{W}^{(1)}_{k,i}$, we get $n(u)=n(u-v)=n(u+v)=k-1$, $n(v)=i-1$ and $u(0)>|v(0)|$;
  \item [$(4).$] for any $i=1,\cdots,k-1$, $(-\infty,\widetilde{\beta}_{k,i})\subset\mbox{Proj}_\beta(\mathcal{W}^{(1)}_{k,i})\subset(-\infty,1)$.  And for any $(\beta,u,v)\in\mathcal{W}^{(1)}_{k,i}$, we get $n(u)=n(u-v)=n(u+v)=k-1$, $n(v)=i-1$ and $u(0)>|v(0)|$;
  \item [$(5).$] the branch $\mathcal{W}^{(1)}_{k,k}\subset\mathcal{O}^1_k$ contains a circle of solutions
  \begin{align}
  \Big\{s_\theta=(1,\cos\theta \, w_k,\sin\theta\,w_k) \Big|\theta\in(0,2\pi)\Big\}.\nonumber
  \end{align}
  Moreover, $\overline{\mathcal{U}^{(l)}_{k,k}} \cap\mathcal{ST}^1_k=\overline{\mathcal{W}^{(1)}_{k,k}}\cap\mathcal{ST}^1_k =\{(1,w_k,0)\}$ for $l=1,2,3,4$.
\end{itemize}
\end{theorem}
A schematic diagram of bifurcations emanating from $\mathcal{T}^1_k$ and $\mathcal{ST}_k^1$ is in Figure \ref{fig:drawing-1-1}.
\begin{remark}\label{r:Wbifurcationparameters}
Analogue results hold for the bifurcating branches emanating out of $\mathcal{ST}_k^2$, $\mathcal{ST}_k^3$ and $\mathcal{ST}_k^4$. We denote the corresponding branches by $\mathcal{W}^{(2)}_{k,i}$, $\mathcal{W}^{(3)}_{k,i}$ and $\mathcal{W}^{(4)}_{k,i}$, respectively. It is easy to verify that the corresponding bifurcation parameters are also $\widetilde{\beta}_{k,i}$.
\end{remark}

As a complement to the above results, we have the following nonexistence result. The novelty of the result is to point out that not all problems with prescribed number of nodes have a solution.
\begin{theorem}\label{t:nonexistencebeta=3}
For any $P,Q\in\mathbb{N}$ with $P\neq Q$, there exists a constant $\varepsilon_0=\varepsilon_0(P,Q)>0$ such that Problem (\ref{e:001}) admits no solution $(u,v)$ such that $n(u)=P$ and $n(v)=Q$ if $\beta\in[3-\varepsilon_0,3+\varepsilon_0]$.
\end{theorem}

\begin{remark}
It is worth to be pointed out that the existence of solutions with \emph{arbitrarily many nodes} holds for Problem (\ref{e:001}) with large $\beta$ (see \cite{MaiaMontefuscoPellacci2008}) or small $\beta$ (see \cite{LiWang2021,LiuWang2019}). Theorem \ref{t:nonexistencebeta=3} refines the above results from the point of view on nonexistence.
\end{remark}

\subsection{An application of Theorem \ref{t:bifurcation}}\label{sub:comparison}
To our surprise, Theorem \ref{t:bifurcation} on bifurcations on the two system can be used to compare different nodal solutions of the scalar field equation
\begin{equation}\label{e:002}
    \left\{
   \begin{array}{lr}
     -{\Delta}u+u= u^3\mbox{ in }B_1 ,\\
     u(0)>0\mbox{ and }n(u)=k-1,\\
     u\in H_{0,r}^1(B_1).
   \end{array}
   \right.
\end{equation}
\begin{theorem}\label{t:comparison}
Let $w_k$ be the unique solution to Problem (\ref{e:002}) with $w(0)>0$ and $n(w)=k-1$.
Then for any $k,i\in\mathbb{N}_+$ with $k>i$, $n(w_k - w_i)=n(w_k + w_i)=k-1$.
\end{theorem}
Since the proof is short and follows from our main result immediately we give it here.
\noindent{\bf Proof.}
This is an immediate consequence of the facts that
\begin{itemize}
  \item for any $i=1,\cdots,k-1$, $0\in\mbox{Proj}_\beta(\mathcal{W}^{(1)}_{k,i})$;
  \item for any $i,k\in\mathbb{N}_+$, the solution of $(u,v)$ to Problem (\ref{e:001}) with $\beta=0$ is unique.
\end{itemize}
\begin{flushright}
$\Box$
\end{flushright}

\subsection{Main results: Liouville type theorems}
To study the bifurcations of the synchronized branches and the semi-trivial branches of the coupled systems, we need to understand the scalar field counterparts in detail. The following result provides the Morse index information which is crucial in the bifurcation analysis.
\begin{proposition}\label{p:MorseIndex}
For any $k\in\mathbb{N}_+$, the solution $w_k\in H_{0,r}^1(B_1)$ to the problem
\begin{equation}
    \left\{
   \begin{array}{lr}
     -{\Delta}u+u=u^3\mbox{ in }B_1 ,\nonumber\\
     u(0)>0\mbox{ and }n(u)=k-1,\nonumber\\
     u\in H_{0,r}^1(B_1)\nonumber
   \end{array}
   \right.
\end{equation}
is unique and non-degenerate and has Morse index $k$ in $H_{0,r}^1(B_1)$. Here, by the word "Morse index of $u$ in $H_{0,r}^1(B_1)$", we mean the dimension of the negative space of $D^2 I(u)$ and $I(u)=\frac{1}{2}\int_{B_1}|\nabla u|^2+u^2-\frac{1}{4}\int_{B_1}|u|^4$ in the space $H_{0,r}^1(B_1)$.
\end{proposition}

We denote the Morse index of $u$ by $m(u)$.
The existence and uniqueness are proved in \cite{WillemBook1996,Tanaka2016}.
The result concerning Morse index is proved in \cite{HarrabiRebhi2011} by an ODE method. In this paper, we will prove Proposition \ref{p:MorseIndex} via a perturbation argument which is given in the Appendix \ref{APPENDIX}.

Next we need new Liouville-type theorems concerning with nodal solutions. Our results focus on the following three homogeneous problems:
\begin{equation}\label{e:Liouville1}
    \left\{
   \begin{array}{lr}
     -{\Delta}u=u^3+\beta uv^2 \,\,\,\,\,\,\,\,  \mbox{in}\ \mathbb{R}^N ,\\
     -{\Delta}v=v^3+\beta u^2v \,\,\,\,\,\,\,\,  \mbox{in}\ \mathbb{R}^N ,\\
     u,v\in C_b^2(\mathbb{R}^N)\mbox{ and }u,v\mbox{ are radial}
   \end{array}
   \right.
\end{equation}
\begin{equation}\label{e:Liouville2}
    \left\{
   \begin{array}{lr}
     -u''=u^3+\beta uv^2 \,\,\,\,\,\,\,\,  \mbox{in}\ \mathbb{R} ,\\
     -v''=v^3+\beta u^2v \,\,\,\,\,\,\,\,  \mbox{in}\ \mathbb{R} ,\\
     u,v\in C_b^2(\mathbb{R})
   \end{array}
   \right.
\end{equation}
and
\begin{equation}\label{e:Liouville3}
    \left\{
   \begin{array}{lr}
     -u''=u^3+\beta uv^2 \,\,\,\,\,\,\,\,  \mbox{in}\ [0,\infty) ,\\
     -v''=v^3+\beta u^2v \,\,\,\,\,\,\,\,  \mbox{in}\ [0,\infty) ,\\
     u,v\in C_b^2([0,\infty))\mbox{ with }u'(0)=v'(0)=0.
   \end{array}
   \right.
\end{equation}
Here, $N=2,3$.
\begin{theorem}\label{t:Liouville}
If $\beta\neq-1$ and if $n(u+v),n(u-v)<\infty$, then Problem (\ref{e:Liouville1}) admits no nontrivial solutions. If $\beta=-1$ and if $n(u+v),n(u-v)<\infty$, then any solution $(u,v)$ of Problem (\ref{e:Liouville1}) are of the form $(c,-c)$ for some constant $c\in\mathbb{R}$.
\end{theorem}
\begin{theorem}\label{t:Liouville2}
If $\beta\neq-1$ and if $n(u+v),n(u-v)<\infty$, then Problems (\ref{e:Liouville2}) and (\ref{e:Liouville3}) admit no non-trivial solutions. If $\beta=-1$ and if $n(u+v),n(u-v)<\infty$, then any solution $(u,v)$ of Problems (\ref{e:Liouville2}) and (\ref{e:Liouville3}) are of the form $(c,-c)$ for some constant $c\in\mathbb{R}$.
\end{theorem}
The proofs of Theorem \ref{t:Liouville} and of Theorem \ref{t:Liouville2} are similar, we only prove Theorem \ref{t:Liouville}.

\subsection{Organization of this paper}
We start in Section 2 to prove our main results Theorems \ref{t:bifurcation} and \ref{t:bifurcation2} which make use of Theorems \ref{t:Liouville} and \ref{t:Liouville2} whose proofs are postponed to Section 4. Theorem \ref{t:nonexistencebeta=3} is proved in Section 3. In Section \ref{section:miscellanies}, we give some additional results on the asymptotics of the global bifurcations. Proposition \ref{p:MorseIndex} is proved in Appendix \ref{APPENDIX}.

A detailed outline for our approach to prove the global bifurcation results in Section \ref{section:bifurcation} is divided into the following steps:
\begin{itemize}
  \item [$(1).$] Compute the bifurcation parameters via Morse indices and determine whether they create global bifurcations (cf. Subsection 2.1-2.3);
  \item [$(2).$] For any $i\in\mathbb{N}_+\backslash\{k\}$, prove that $\mathcal{U}^{(1)}_{k,i}$ ($\mathcal{W}^{(1)}_{k,i}$) does not intersect with any other branches. This is done by proving that the numbers of nodes of $u$, $v$, $u+v$ and $u-v$ are constant on each of the branches (cf. Subsection \ref{sub:assertion23});
  \item [$(3).$] Prove that $\mathcal{U}^{(1)}_{k,i}$ ($\mathcal{W}^{(1)}_{k,i}$) is bounded on each compact interval of $\beta$-axis by using Liouville-type theorems, and determine the projection range of $\mbox{Proj}_\beta(\mathcal{U}^{(1)}_{k,i})$ ($\mbox{Proj}_\beta(\mathcal{W}^{(1)}_{k,i})$) (cf. Subsection \ref{sub:assertion23});
  \item [$(4).$] Verify that $\mathcal{U}^{(1)}_{k,i}$ ($\mathcal{W}^{(1)}_{k,i}$, respectively) intersects with the main branches (cf. Subsection \ref{subsection:(5)}).
\end{itemize}

\section{Proof of Theorems \ref{t:bifurcation} and \ref{t:bifurcation2}}\label{section:bifurcation}

In this part, we follow an idea in \cite{BartschDancerWang2010} to prove Theorem \ref{t:bifurcation}.
First, let us define some notions. $m(\beta)$ denotes the Morse index of $J_\beta$ at $(u_\beta,v_\beta)$.

\subsection{A weighted eigenvalue problem}
As in \cite{BartschDancerWang2010}, our main approach to determine the bifurcation parameters is the computation of Morse index. To this end, we need to study the following weighted eigenvalue problem.
\begin{equation}\label{e:003}
    \left\{
   \begin{array}{lr}
     -{\Delta}u+u=\lambda w_k^2u\mbox{ in }B_1 ,\\
     u\in H_{0,r}^1(B_1)
   \end{array}
   \right.
\end{equation}
Here, the $w_k$ is the solution to Problem (\ref{e:002}). The following claims hold evidently.
\begin{claim}\label{c:uniqueness}
The solution $w_k$ to Problem (\ref{e:002}) is unique and non-degenerate.
\end{claim}
This is due to \cite{Tanaka2016,LiMiyagaki}.
\begin{claim}\label{c:discrete}
Problem (\ref{e:003}) admits only discrete spectrum.
\end{claim}
Claim \ref{c:discrete} holds evidently.
\begin{claim}\label{claim:3}
$\lambda = 1$ is an eigenvalue of Problem (\ref{e:003}) with corresponding eigenfunction $w_k$. $\lambda =3$ is not an eigenvalue of Problem (\ref{e:003}).
\end{claim}
The first part of Claim \ref{claim:3} is obvious. And second part is due to the non-degeneracy of $w_k$.

\begin{definition}
Denote the set $\sigma_k=\{\lambda\in\mathbb{R}|\lambda\mbox{ is an eigenvalue of Problem }(\ref{e:003})\}$.
\end{definition}

The above claims imply the following lemma.
\begin{lemma}\label{l:eigenvalues}
The following hold for Problem (\ref{e:003}).
\begin{itemize}
  \item [$(1).$] The eigenvalues $\{\lambda_{k,i}\}_{i=1}^\infty$ of Problem (\ref{e:003}) satisfy that
  $$0<\lambda_{k,1}<\cdots<\lambda_{k,k}=1<3<\lambda_{k,k+1}<\lambda_{k,k+2}<\cdots;$$
  \item [$(2).$] For any $i\in\mathbb{N}$, the corresponding eigenfunction $\phi_{k,i}$ of $\lambda_{k,i}$ changes its sign exactly $i-1$ times. Especially, $\phi_{k,1}>0$ in $B_1$ and $\phi_{k,k}=w_k$.
\end{itemize}
\end{lemma}

\noindent{\bf Proof.}
Firstly, it is obvious that all of the eigenspace are one dimensional. This is due to the radial symmetry of
Problem (\ref{e:003}).

By Proposition \ref{p:MorseIndex}, $\#\{\lambda<3|\lambda\in\sigma_k\}=k$. On the other hand, \cite[Theorem XIII.7.53 and Corollary 7.56]{DunfordSchwartzBook1988} implies that $\phi_{k,i}$ changes its sign exactly $i-1$ times. Therefore, $w_k$ is the corresponding eigenfunction of $\lambda_{k,k}$, since its nodal number $n(w_k)=k-1$. And we also have that $\lambda_{k,k}=1$. Therefore, $0<\lambda_{1,k}<\cdots<\lambda_{k,k}=1$ and $\#\{\lambda\leq1|\lambda\in\sigma_k\}=k$, which implies that $(1,3)\cap\sigma_k=\emptyset$.

In order to complete the proof, we only need to show that $\lambda_{k,k+1}>3$. This is valid since $3\notin\sigma_k$, cf. Claim \ref{claim:3}.

\begin{flushright}
$\Box$
\end{flushright}

\subsection{Determining the bifurcation parameters}
According to \cite[Section 8.9]{MawhinWillem1989}, the bifurcation can only occur where the linearization of the operator possesses a nontrivial kernel. Therefore, the first step is to determine the parameters $\beta$ that ensure the non-triviality of the kernel. To this end, we consider the linearzation of Problem (\ref{e:001}) at $(\beta,u_\beta,v_\beta)\in\mathcal{T}^1_k$.
\begin{equation}\label{e:linear1}
    \left\{
   \begin{array}{lr}
     -{\Delta}\phi+\phi=3u_\beta^2\phi+\beta v_\beta^2\phi+2\beta u_\beta v_\beta\psi\mbox{ in }B_1 ,\\
     -{\Delta}\psi+\psi=2\beta u_\beta v_\beta\phi+3v_\beta^2\psi+\beta u_\beta^2\psi\mbox{ in }B_1,\\
     \phi,\psi\in H_{0,r}^1(B_1),
   \end{array}
   \right.
\end{equation}
i.e.,
\begin{equation}\label{e:linear}
    \left\{
   \begin{array}{lr}
     -{\Delta}\phi+\phi=w_k^2\Big(\frac{\beta+3}{1+\beta}\phi+\frac{2\beta}{1+\beta}\psi\Big)\mbox{ in }B_1 ,\\
     -{\Delta}\psi+\psi=w_k^2\Big(\frac{2\beta}{1+\beta}\phi+\frac{\beta+3}{1+\beta}\psi\Big)\mbox{ in }B_1 ,\\
     \phi,\psi\in H_{0,r}^1(B_1).
   \end{array}
   \right.
\end{equation}
This gives
\begin{equation}
    \left\{
   \begin{array}{lr}
     -{\Delta}(\phi+\psi)+(\phi+\psi)=3w_k^2(\phi+\psi)\mbox{ in }B_1 ,\nonumber\\
     -{\Delta}(\phi-\psi)+(\phi-\psi)=\frac{3-\beta}{1+\beta}w_k^2 (\phi-\psi)\mbox{ in }B_1 ,\nonumber\\
     \phi,\psi\in H_{0,r}^1(B_1).\nonumber
   \end{array}
   \right.
\end{equation}
Since $w_k$ is non-degenerate, $\phi+\psi\equiv0$. Then, the above equations admits a nontrivial solution if and only if
\begin{align}\label{InclusionEigen}
\frac{3-\beta}{1+\beta}\in\sigma_k.
\end{align}
Here, $\sigma_k$ is the set of spectrum of Problem (\ref{e:003}).
Inclusion (\ref{InclusionEigen}) is equivalent to
\begin{align}\label{BifurcationParameters}
\beta=\beta_{k,i}:=\frac{4}{\lambda_{k,i}+1}-1
\end{align}
with $\lambda_{k,i}$ is the $i$-th eigenvalue of Problem (\ref{e:003}).
The solution set of Problem (\ref{e:linear}) with $\beta=\beta_{k,i}$ is
\begin{align}\label{Vl}
V_i:=\{c\cdot(\phi_{k,i},-\phi_{k,i})|c\in\mathbb{R},\,i\in\mathbb{N}\mbox{ and }\phi_{k,i}\mbox{ is the }i\mbox{-th eigenfunction of Probelm (\ref{e:003})}\}.
\end{align}
In such a way, we prove the following claim.
\begin{claim}
$\{\beta_{k,i}|\beta_{k,i}\mbox{ is defined as (\ref{BifurcationParameters})}\}$ is the set of possible bifurcation parameters of $\mathcal{T}^1_k$.
\end{claim}

Moreover, the following result is obvious.
\begin{claim}
It holds that
\begin{align}
-1<\cdots<\beta_{k,k+2}<\beta_{k,k+1}<0<1 =\beta_{k,k}<\beta_{k,k-1}<\cdots<\beta_{k,1}<3.\nonumber
\end{align}
\end{claim}

\subsection{Proof of Assertion $(1)$ of Theorem \ref{t:bifurcation} and of Assertion $(1)$ of Theorem \ref{t:bifurcation2}}
In this part, we show the bifurcations of $\mathcal{T}^1_k$ occur at each $\beta=\beta_{k,i}$.
To this end, we need to show the jumping of Morse indices.
\begin{align}
H_\beta[\phi,\psi]=\int|\nabla \phi|^2+|\phi|^2+\int|\nabla\psi|^2+|\psi|^2-\int w_k^2 \Big[\frac{\beta+3}{\beta+1}(\phi^2+\psi^2)+\frac{4\beta}{\beta+\mu}\phi\psi\Big].
\end{align}
The derivative of $H_\beta$ in $\beta$ is that
\begin{align}
\frac{d}{d\beta}H_\beta[\phi,\psi]=\frac{2}{(\beta+1)^2}\int(\phi^2+\psi^2-2\phi\psi).\nonumber
\end{align}

Notice that $H_{0,r}^1(B_1)\times H_{0,r}^1(B_1)=\mbox{span}\big\{(\phi_{k,i},\phi_{k,i}),(\phi_{k,i},-\phi_{k,i})\big| i\in\mathbb{N}_+\big\}$. To compute the Morse index $m(\beta)$, it is sufficient to check the sign of $H_\beta[\phi_{k,i},\phi_{k,i}]$ and $H_\beta[\phi_{k,i},-\phi_{k,i}]$ for any $i\in\mathbb{N}_+$. Here, $\phi_{k,i}$ is the $i$-th eigenfunction of Problem (\ref{e:003}).
First, for any $(\phi_{k,i},\phi_{k,i})$ and any $\beta\in(-1,3)$,
\begin{align}
H_\beta[\phi_{k,i},\phi_{k,i}]=2\|\phi_{k,i}\|^2-6\int\phi_{k,i}^2 w_k^2.\nonumber
\end{align}
Then, $H_\beta[\phi_{k,i},\phi_{k,i}]<0$ for $i=1,\cdots,k$ and $H_\beta[\phi_{k,i},\phi_{k,i}]>0$ for $i=k+1,\cdots$. For $(\phi_{k,i},-\phi_{k,i})$,
\begin{align}
H_\beta[\phi_{k,i},-\phi_{k,i}]=2\|\phi_{k,i}\|^2-2\cdot\frac{3-\beta}{\beta+1}\int\phi_{k,i}^2 w_k^2.\nonumber
\end{align}
Recall that we define $\beta_{k,i}=\frac{4}{\lambda_{k,i} +1}-1$ with $\lambda_{k,i}\in\sigma_k$ is the $i$-th eigenvalue of Problem (\ref{e:003}). If $\beta\in[\beta_{k,i+1},\beta_{k,i})$, we have $H_{\beta}[\phi_{k,j},-\phi_{k,j}]<0$ for $j=1,\cdots,i-1$. $H_{\beta}[\phi_{k,i},-\phi_{k,i}]\geq0$ and the equality holds if and only if $\beta=\beta_{k,i+1}$. In this case, $(\phi_{k,i},-\phi_{k,i})\in\mbox{Ker}(H_{\beta})=\mbox{Ker}(H_{\beta_{k,i+1}})$.
On the other hand, $H_{\beta}[\phi_{k,j},-\phi_{k,j}]>0$ for any $j=i+2,\cdots$.

The above computations also give that
\begin{claim}
The kernel space of $H_{\beta_{k,i}}$ is $V_{k,i}$ defined in (\ref{Vl}) for $i\in\mathbb{N}_+$.
\end{claim}
Then, we have the following corollary.
\begin{claim}
For any $i\in\mathbb{N}_+$, there is a constant $\varepsilon_{0,i}>0$ such that for any $\varepsilon\in(0,\varepsilon_{0,i})$, $H_{\beta_{k,i}-\varepsilon}$ is negative definite on $V_{k,i}$.
\end{claim}

\noindent{\bf Proof.}
To verify this, one only needs to notice that $\frac{d}{d\beta}H_\beta|_{\beta=\beta_{k,i}}>0$ on $V_{k,i}$.
\begin{flushright}
$\Box$
\end{flushright}

In a word, we proved the following claim.
\begin{claim}
For any $\beta\in[\beta_{k,i+1},\beta_{k,i})$,
$m(\beta)=k+i$. And $\dim(\mbox{Ker}(H_\beta))=1$ if $\beta=\beta_{k,i}$ for some $i\in\mathbb{N}_+$. Otherwise, $\dim(\mbox{Ker}(H_\beta))=0$.
\end{claim}

An immediate consequence is that
\begin{corollary}
For any $i\in\mathbb{N}_+$, $\lim_{\varepsilon\to0+}\Big(m(\beta_{k,i} +\varepsilon)-m(\beta_{k,i}-\varepsilon)\Big)=-1$.
\end{corollary}
Therefore, we have the following bifurcation results.
\begin{corollary}\label{coro:globalbranch}
For any $i\in\mathbb{N}_+$, $\beta=\beta_{k,i}$ is a bifurcation parameter. Especially, there exists a connected subset $\mathcal{U}^{(1)}_{k,i}$ of $\mathcal{S}^1_k\cup\{(\beta_{k,i},u_{\beta_{k,i}},v_{\beta_{k,i}})\}$ either
\begin{itemize}
  \item [$(a).$] $\mathcal{U}^{(1)}_{k,i}$ is unbounded in $\mathbb{R}\times H_{0,r}^1(B_1)\times H_{0,r}^1(B_1)$, or
  \item [$(b).$] there exist a $(\beta_{k,i'},u_{\beta_{k,i'}} ,v_{\beta_{k,i'}})\in\mathcal{T}^1_k$ with $\frac{3-\beta_{k,i'}}{1+\beta_{k,i'}}\in\sigma_k$ and $i\neq i'$ such that $(\beta_{k,i'},u_{\beta_{k,i'}},v_{\beta_{k,i'}})\in\overline{\mathcal{U}^{(1)}_{k,i}}$.
\end{itemize}
\end{corollary}
This is due to \cite[Theorem 1.3]{Rabinowitz1971}.
In the following, we will prove that for any $(\beta,u,v)\in\mathcal{U}^{(1)}_{k,i}$, $n(u)=n(v)=k$.
Moreover, notice that $\beta=\beta_{k,i}$ is a bifurcation parameter of a simple eigenvalue, we use \cite{CrandallRabinowitz1971} to get the following result.

\begin{lemma}\label{l:1parameter}
For any $i\in\mathbb{N}_+$, there is a constant $r_{k,i}>0$ such that $\mathcal{U}^{(1)}_{k,i}\cap\mathcal{B}_{r_{k,i}}(\beta_{k,i},u_{\beta_{k,i}},v_{\beta_{k,i}})$ can be written as
\begin{align}
(u_{\beta_{k,i}},v_{\beta_{k,i}}) +(\beta-\beta_{k,i})(\phi_{k,i},-\phi_{k,i})+o(\beta-\beta_{k,i}).\nonumber
\end{align}
Here, $\mathcal{B}_{r_{k,i}}(\beta_{k,i}, u_{\beta_{k,i}},v_{\beta_{k,i}})$ is the ball in $\mathbb{R}\times H_{0,r}^1(B_1)\times H_{0,r}^1(B_1)$, centered at $(\beta_{k,i},u_{\beta_{k,i}},v_{\beta_{k,i}})$ with radii $r_{k,i}$.
\end{lemma}
Applying the bootstrap method, cf. \cite{BrezisKato1979}, we have that $o(\beta-\beta_{k,i})$ is $C^{2,\alpha}$-infinitesimal for some $\alpha\in(0,1)$. Therefore, an immediate consequence is following.
\begin{corollary}\label{coro:nodes}
For any $i\in\mathbb{N}_+$, there exists a constant $r_{k,i}$ such that for any $(\beta,u,v)\in \mathcal{U}^{(1)}_{k,i}\cap\mathcal{B}_{r_{k,i}}(\beta_{k,i},u_{\beta_{k,i}} ,v_{\beta_{k,i}})$, $n(u)=n(u_{\beta_{k,i}})=k-1$, $n(v)=v(v_{\beta_{k,i}})=k-1$, $n(u-v)=n(\phi_{k,i})=i-1$, $n(u+v)=n(u_{\beta_{k,i}}+v_{\beta_{k,i}})=k-1$ and  $u(0),v(0)>0$.
\end{corollary}

By a similar approach, we can compute the bifurcation parameters of the semi-trivial branch $\mathcal{ST}_k^1$. The linearization of Problem (\ref{e:001}) at $(w_k,0)$ is that
\begin{equation}\label{e:linear12}
    \left\{
   \begin{array}{lr}
     -{\Delta}\phi+\phi=3w_k^2\phi\mbox{ in }B_1 ,\\
     -{\Delta}\psi+\psi=\beta w_k^2\psi\mbox{ in }B_1,\\
     \phi,\psi\in H_{0,r}^1(B_1).
   \end{array}
   \right.
\end{equation}
Here, $\phi=0$ and $\psi\neq0$ if and only if $\beta=\widetilde{\beta}_{k,i}:=\lambda_{k,i}$ for some $\lambda_{k,i}\in\sigma_k$. Moreover, if $\beta=\widetilde{\beta}_{k,i}$, $\psi=\phi_{k,i}$, the $i$-th eigenfunction of Problem (\ref{e:003}).
By a similar approach, the following results are obtained.
\begin{lemma}\label{l:Semiglobalbranch}
For any $i\in\mathbb{N}_+$, $\beta=\widetilde{\beta}_{k,i}$ is a bifurcation parameter of $\mathcal{ST}_k^1$. Especially, there exists a connected subset $\mathcal{W}^{(1)}_{k,i}$ of $\mathcal{O}^{1}_k\cup\{(\widetilde{\beta}_{k,i},w_k,0)\}$ either
\begin{itemize}
  \item [$(a).$] $\mathcal{W}^{(1)}_{k,i}$ is unbounded in $\mathbb{R}\times H_{0,r}^1(B_1)\times H_{0,r}^1(B_1)$, or
  \item [$(b).$] there exist a $(\beta_{k,i'},w_k,0)\in\mathcal{ST}^1_k$ with $\widetilde{\beta}_{k,i'}=\lambda_{k,i'}$ and $\lambda_{k,i'}\in\sigma_k$ and $i\neq i'$ such that $(\widetilde{\beta}_{k,i'},w_k,0)\in\overline{\mathcal{W}^{(1)}_{k,i}}$.
\end{itemize}
\end{lemma}

\begin{lemma}\label{l:1-parameterSemiTrivial}
For any $i\in\mathbb{N}_+$, there is a constant $\widetilde{r}_{k,i}>0$ such that any $(\beta,u,v)\in\mathcal{W}^{(1)}_{k,i}\cap\mathcal{B}_{\widetilde{r}_{k,i}} (\widetilde{\beta}_{k,i},w_k,0)$ can be written as
\begin{align}
(w_k,0)+(\beta-\widetilde{\beta}_{k,i})(0,\phi_{k,i})+o(\beta-\widetilde{\beta}_{k,i}).\nonumber
\end{align}
Here, $\mathcal{B}_{\widetilde{r}_{k,i}} (\widetilde{\beta}_{k,i},w_k,0)$ is the ball in $\mathbb{R}\times H_{0,r}^1(B_1)\times H_{0,r}^1(B_1)$, centered at $(\widetilde{\beta}_{k,i},w_k,0)$ with radii $\widetilde{r}_{k,i}$. Notice that $o(\beta-\widetilde{\beta}_{k,i})$ is the infinitesimal in $C^{2,\alpha}$-norm for some $\alpha\in(0,1)$.
\end{lemma}
An immediate corollary is that
\begin{corollary}\label{l:SemiTrivialNodes}
For any $i\in\mathbb{N}_+$, there is a constant $\widetilde{r}_{k,i}>0$ such that for any $(\beta,u,v)\in\mathcal{W}^{(1)}_{k,i}\cap \mathcal{B}_{\widetilde{r}_{k,i}}(\widetilde{\beta}_{k,i},w_k,0)$, $n(u)=n(u+v)=n(u-v)=k-1$, $n(v)=i-1$ and $u(0)>|v(0)|$.
\end{corollary}
\begin{remark}\label{r:SemiBranchParameters}
Analogues to Lemmas \ref{l:Semiglobalbranch} and \ref{l:1-parameterSemiTrivial} for the branches $\mathcal{ST}_k^2$, $\mathcal{ST}_k^2$ and $\mathcal{ST}_k^2$ hold with the same bifurcation parameters $\widetilde{\beta}_{k,i}$ and radii $\widetilde{r}_{k,i}$.
\end{remark}
\noindent{\bf Proof of Assertion $(1)$ of Theorem \ref{t:bifurcation} and of Assertion $(1)$ of Theorem \ref{t:bifurcation2}.}
This follows immediately from the above computations.
\begin{flushright}
$\Box$
\end{flushright}

\subsection{Mappings from $\mathcal{T}^{1}_k$ to $\mathcal{ST}_k^l$ for $l=1,2,3,4$}
In this subsection, we define mappings from the synchronized branch $\mathcal{T}^1_k$ to the semi-trivial branches $\mathcal{ST}^l_k$ for $l=1,2,3,4$. In such a way, we transform the result of Theorem \ref{t:bifurcation} to the result of Theorem \ref{t:bifurcation2}. Consider
\begin{align}
u'=\frac{\sqrt{1+\beta}}{2}(u+v)\mbox{ and }v'=\frac{\sqrt{1+\beta}}{2}(u-v).\nonumber
\end{align}
It is easy to check that
\begin{equation}\label{e:comparison}
    \left\{
   \begin{array}{lr}
     -{\Delta}u'+u'= (u')^3+\frac{3-\beta}{1+\beta} u'(v')^2\mbox{ in }B_1 ,\\
     -{\Delta}v'+v'= (v')^3+\frac{3-\beta}{1+\beta} (u')^2v'\mbox{ in }B_1 ,\\
     u',v'\in H_{0,r}^1(B_1).
   \end{array}
   \right.
\end{equation}
Moreover, the function $f(\beta)=\frac{3-\beta}{1+\beta}$ defines a homomorphism from $(-1,+\infty)$ to itself. Especially, $\lim_{\beta\downarrow-1}f(\beta)=+\infty$, $f(0)=3$, $f(1)=1$, $f(3)=0$ and $\lim_{\beta\uparrow+\infty}f(\beta)=-1$.

Define the map
\begin{align}\label{map:T}
T_{1}:(-1,+\infty)\times H_{0,r}^1(B_1)\times H_{0,r}^1(B_1)\to  (-1,+\infty)\times H_{0,r}^1(B_1)\times H_{0,r}^1(B_1)
\end{align}
as
$T_{1}(\beta,u,v)=\Big(\frac{3-\beta}{1+\beta},\frac{\sqrt{1+\beta}}{2}(u+v) ,\frac{\sqrt{1+\beta}}{2}(u-v)\Big)$.
An evident result is that
\begin{lemma}
\begin{itemize}
  \item [$(1)$] 
  $T_{1}(\{(\beta,u,v)\in\mathcal{S}|\beta\in(-1,1)\})= \{(\beta,u,v)\in\mathcal{S}|\beta>1\}$. Here, the set $\mathcal{S}$ is defined as in (\ref{SolutionSet});
  \item [$(2)$] $T_{1}(\mathcal{T}^1_k)=\mathcal{ST}^1_k|_{\beta>-1}$.
\end{itemize}
\end{lemma}
Similarly, for $l=2,3,4$, we define the map
\begin{align}\label{map:T'}
T_{l}:(-1,+\infty)\times H_{0,r}^1(B_1)\times H_{0,r}^1(B_1)\to  (-1,+\infty)\times H_{0,r}^1(B_1)\times H_{0,r}^1(B_1)
\end{align}
as
\begin{align}
T_{2}(\beta,u,v)&=\Big(\frac{3-\beta}{1+\beta},\frac{\sqrt{1+\beta}}{2}(u-v) ,\frac{\sqrt{1+\beta}}{2}(u+v)\Big)\nonumber\\
T_{3}(\beta,u,v)&=\Big(\frac{3-\beta}{1+\beta},-\frac{\sqrt{1+\beta}}{2}(u+v) ,\frac{\sqrt{1+\beta}}{2}(u-v)\Big)\nonumber
\end{align}
and
\begin{align}
T_{4}(\beta,u,v)=\Big(\frac{3-\beta}{1+\beta},\frac{\sqrt{1+\beta}}{2}(u-v) ,-\frac{\sqrt{1+\beta}}{2}(u+v)\Big).\nonumber
\end{align}
\begin{lemma}
For $l=2,3,4$, we get
\begin{itemize}
  \item [$(1)$] $T_{l}(\{(\beta,u,v)\in\mathcal{S}|\beta\in(-1,1)\})= \{(\beta,u,v)\in\mathcal{S}|\beta>1\}$. Here, the set $\mathcal{S}$ is defined as in (\ref{SolutionSet});
  \item [$(2)$] $T_{l}(\mathcal{T}_k^1)=\mathcal{ST}^l_k|_{\beta>-1}$.
\end{itemize}
\end{lemma}

\subsection{Proof of Assertions $(2-3)$ of Theorem \ref{t:bifurcation}}\label{sub:assertion23}

First, we need to check that

\begin{lemma}\label{l:Sign1}
Assume $i=k+1,k+2,\cdots$.
For any $(\beta,u,v)\in\mathcal{U}^{(1)}_{k,i}$, we have that $n(u)=n(v)=n(u+v)=k-1$, $n(u-v)=i-1$, $u(0)>0$ and $v(0)>0$.
\end{lemma}
\begin{remark}
In Subsection \ref{subsection:(5)}, we will see that the assumption $i\neq k$ is necessary.
\end{remark}

Before proving Lemma \ref{l:Sign1}, we note the following topological result.

\begin{claim}\label{c:FinitelyManyComponents}
$\mathcal{U}^{(1)}_{k,i}\cup\{(\beta_{k,i},u_{\beta_{k,i}},v_{\beta_{k,i}})\}$, $\{(\beta,u,v)\in \mathcal{U}^{(1)}_{k,i}\cup\{(\beta_{k,i}, u_{\beta_{k,i}},v_{\beta_{k,i}})\}|\beta\geq0\}$ and $\{(\beta,u,v)\in \mathcal{U}^{(1)}_{k,i}\cup\{(\beta_{k,i},u_{\beta_{k,i}},v_{\beta_{k,i}})\}|\beta\leq0\}$
contain at most countable connected components.
\end{claim}

\noindent{\bf Proof.}
We only need to check the claim holds for $\mathcal{U}^{(1)}_{k,i}\cup\{(\beta_{k,i},u_{\beta_{k,i}} ,v_{\beta_{k,i}})\}$. First of all, let us notice that we only need to show that $\mathcal{U}^{(1)}_{k,i}\cup\{(\beta_{k,i},u_{\beta_{k,i}},v_{\beta_{k,i}})\}$ can be represented as a union of countable of compact sets. The claim follows from the fact that $\mathbb{R}\times H_{0,r}^1(B_1)\times H_{0,r}^1(B_1)$ is a Hausdorff space.

Define the set $\mathcal{C}_n:=\overline{\mathcal{B}_n(0,0,0)\backslash \mathcal{B}_{n-1}(0,0,0)}\cap \Big(\mathcal{U}_{k,i}^{(1)}\cup\{(\beta_{k,i},u_{\beta_{k,i}},v_{\beta_{k,i}})\}\Big)$. Then
\begin{itemize}
  \item [$(1).$] $\mathcal{C}_n$ is compact due to the $(PS)$ condition;
  \item [$(2).$] $\cup_{n\geq1}\mathcal{C}_n=\mathcal{U}_{k,i}^{(1)}\cup \{(\beta_{k,i},u_{\beta_{k,i}},v_{\beta_{k,i}})\}$.
\end{itemize}
Claim \ref{c:FinitelyManyComponents} follows immediately.
\begin{flushright}
$\Box$
\end{flushright}

Our proof of Lemma \ref{l:Sign1} is divided into two steps. The first step proves that if Lemma \ref{l:Sign1} fails, $\mathcal{U}^{(1)}_{k,i}$ will intersect with one of the semi-trivial branches  $\mathcal{ST}_{k'}^l$ for $l=1,2,3,4$ and for some $k'\in\mathbb{N}_+$. The second step proves that $\mathcal{U}^{(1)}_{k,i}$ cannot intersect with any semi-trivial branch, which completes the proof.

~~

\noindent{\bf Proof.}
We assume that $i=k+1,k+2,\cdots$.

First of all, it is known by Corollary \ref{coro:nodes} that for any $(\beta,u,v)\in\mathcal{U}^{(1)}_{k,i}\cap\mathcal{B}_{r_{k,i}}(\beta_{k,i} ,u_{\beta_{k,i}},v_{\beta_{k,i}})$ we have $n(u)=n(v)=n(u+v)=k-1$ and $n(u-v)=i-1$. We argue by contradiction and prove that if Lemma \ref{l:Sign1} fails, then we only need to consider the intersection of $\mathcal{U}^{(1)}_{k,i}$ with a semi-trivial branch $\mathcal{ST}_{k'}^l$ for some $l=1,2,3,4$.

~~

\noindent{\bf Step 1. Only need to check the intersection of $\mathcal{U}^{(1)}_{k,i}$ with semi-trivial branches.}

If Lemma \ref{l:Sign1} does not hold, then there is a point $(\beta',u',u')\in\mathcal{U}^{(1)}_{k,i}$ with at least one of the following holds:
\begin{itemize}
  \item [$(1).$] $u'(0)\leq 0$;
  \item [$(2).$] $v'(0)\leq 0$.
  \item [$(3).$] $n(u')\neq k-1$;
  \item [$(4).$] $n(v')\neq k-1$;
  \item [$(5).$] $n(u'+v')\neq k-1$;
  \item [$(6).$] $n(u'-v')\neq i-1$.
\end{itemize}
First, if $(1)$ or $(2)$ occurs, there exists a $(\beta'',u'',v'')\in\mathcal{U}^{(1)}_{k,i}$ such that $u''=0$ or $v''=0$. Indeed, since any $(\beta,u,v)\in\mathcal{U}^{(1)}_{k,i}$ is a radial solution to Problem (\ref{e:001}), $\frac{d}{dr}u(0)=\frac{d}{dr}v(0)=0$. On the other hand, by the connectedness of $\mathcal{U}^{(1)}_{k,i}$, there exists a $(\beta'',u'',v'')\in\mathcal{U}^{(1)}_{k,i}$ such that $u''(0)=0$ or $v''(0)=0$. Then $u''=0$ or $v''=0$ follows immediately. Since $(0,0)$ is isolated, then $(u'',v'')=(0,v'')$ or $(u'',v'')=(u'',0)$. This implies that $\mathcal{U}^{(1)}_{k,i}$ intersects with a semi-trivial branch.

If Cases $(3)$ or $(4)$ occurs, there is a $(\beta'',u'',v'')\in\mathcal{U}^{(1)}_{k,i}$ such that $u''=0$ or $v''=0$.

In Case $(5)$, there is a $(\beta'',u'',v'')\in\mathcal{U}^{(1)}_{k,i}$ such that $u''+v''=0$. In this case, $u''=-v''$. Notice that $v''(0)>0$, then $u''(0)=\frac{d}{dr}u''(0)=0$ or $v''(0)=\frac{d}{dr}v''(0)=0$. Therefore, we again find a $(\beta''',u''',v''')\in\mathcal{U}^{(1)}_{k,i}$ with either $u'''=0$ or $v'''=0$.

In Case $(6)$, there is a $(\beta'',u'',v'')\in\mathcal{U}^{(1)}_{k,i}$ with $u''=v''$. Then, $\mathcal{U}^{(1)}_{k,i}\cap\mathcal{T}^1_{k'}\neq\emptyset$ for some $k'\neq k$ or $\mathcal{U}^{(1)}_{k,i}$ returns to $\mathcal{T}^1_k$. If $\mathcal{U}^{(1)}_{k,i}\cap\mathcal{T}^1_{k'}\neq\emptyset$ for some $k'\neq k$, there is a $(\beta''',u''',v''')\in\mathcal{U}^{(1)}_{k,i}$ with $u'''=0$ or $v'''=0$. Consider the isolation of the trivial solution, we only need to check the case that $(\beta''',u''',v''')\in\mathcal{ST}_{k'''}^l$ for some $l=1,2,3,4$ and $k'''\in\mathbb{N}$.
If $\mathcal{U}^{(1)}_{k,i}$ returns to $\mathcal{T}^1_k$, then we denote
\begin{itemize}
  \item $\mathcal{A}_1=\cup_{j\neq i}\overline{\mathcal{B}_{\frac{r_{k,j}}{4}}(\beta_{k,j}, u_{\beta_{k,j}},v_{\beta_{k,j}})} \cap\mathcal{U}^{(1)}_{k,i}$;
  \item $\mathcal{B} =\overline{\mathcal{U}^{(1)}_{k,i} \cup\{(\beta_{k,i},u_{\beta_{k,i}},v_{\beta_{k,i}})\}\backslash\mathcal{A}_1}$;
  \item $\mathcal{A}_2=\cup_{j\neq i}\overline{\mathcal{B}_{\frac{r_{k,j}}{2}}(\beta_{k,j},u_{\beta_{k,j}} ,v_{\beta_{k,j}})} \cap\mathcal{U}^{(1)}_{k,i}$.
\end{itemize}

Here, $r_{k,j}$'s are the radius in Lemma \ref{l:1parameter}.
Let us check that $\mathcal{B}\cap\mathcal{A}_2\neq\emptyset$. If $\mathcal{B}\cap\mathcal{A}_2=\emptyset$, then we have that $\mathcal{B}\cap\mathcal{A}_1=\emptyset$. Moreover, it is evident that $\mathcal{U}^{(1)}_{k,i}\cup\{(\beta_{k,i},u_{\beta_{k,i}},v_{\beta_{k,i}})\}=\mathcal{B}\cup\mathcal{A}_1$. This concludes that $\mathcal{U}^{(1)}_{k,i}$ is disconnected. This is a contradiction. Therefore, $\mathcal{B}\cap\mathcal{A}_2\neq\emptyset$.

On the other hand, similar to Claim \ref{c:FinitelyManyComponents}, the set $\mathcal{B}$ contains at most countable connected components. Then there is a connected component of $\mathcal{B}$, say $B_{*}$, with $(\beta_{k,i},u_{\beta_{k,i}},v_{\beta_{k,i}})\in B_{*}$ and $B_{*}\cap\mathcal{A}_2\neq\emptyset$. If $(\beta_{k,i},u_{\beta_{k,i}},v_{\beta_{k,i}})\in B_{*}$ with $B_{*}\cap\mathcal{A}_2=\emptyset$, we again obtain the disconnectedness of $\mathcal{U}^{(1)}_{k,i}\cup\{(\beta_{k,i},u_{\beta_{k,i}},v_{\beta_{k,i}})\}$ by proving it can be divided into two disjoint nonempty closed sets, $B_{*}$ and $(\mathcal{B}\backslash B_*)\cup\mathcal{A}_1$.
Due to Corollary \ref{coro:globalbranch}, there exists a $(\beta'',u'',v'')\in B_*\cap\mathcal{A}_2$ such that $n(u'')=n(v'')=n(u''+v'')=k-1$ and $n(u''-v'')\neq i-1$. Then, we can find a $(\beta''',u''',v''')\in B_*$ with $u'''=v'''$.
If $u'''=v'''=0$, this contradicts with the isolation of the trivial solution $(0,0)$. If $u'''=v'''\neq 0$, then $\beta>-1$ due to the explicit formulation of Problem (\ref{e:001}). Due to the uniqueness result (cf. \cite{Tanaka2016}), $(\beta''',u''',v''')\in\mathcal{T}^1_{n(u''')+1}$. Notice that if $u'''(0)\leq 0$ or $v'''(0)\leq0$, we return to Case $(1)$ or Case $(2)$.
Due to the construction of $B_*$, $B_*\backslash\{(\beta_{k,i},u_{\beta_{k,i}},v_{\beta_{k,i}})\}\cap\mathcal{T}^1_k=\emptyset$.
Then, $n(u''')=n(v''')\neq k$. This returns to the first alternative. Again, we prove that $\mathcal{U}^{(1)}_{k,i}$ intersects with a semi-trivial branch.

~~

\noindent{\bf Step 2. $\mathcal{U}^{(1)}_{k,i}$ cannot intersect with any semi-trivial branch.}

The above argument conclude that $\mathcal{U}^{(1)}_{k,i}\cap\mathcal{ST}^l_{k'}\neq\emptyset$ for some $l=1,2,3,4$ and some $k'\in\mathbb{N}$. Therefore, in order to induce a contradiction, we only need to consider the intersection of $\mathcal{U}^{(1)}_{k,i}$ with some semi-trivial branch $\mathcal{ST}_{k'}^l$.

For any $i\geq k+1$ and any $k'\in\mathbb{N}$, all the bifurcation parameters of $\mathcal{ST}_{k'}^1$, $\mathcal{ST}_{k'}^2$, $\mathcal{ST}_{k'}^3$ and $\mathcal{ST}_{k'}^4$ are on the semiaxis $\{\beta>0\}$. On the other hand, the bifurcation parameter $\beta_{k,i}$ of $\mathcal{U}^{(1)}_{k,i}$ is negative. Then, $0\in\mbox{Proj}_\beta(\mathcal{U}^{(1)}_{k,i})$. By Claim \ref{c:FinitelyManyComponents},  $\{(\beta,u,v)\in\mathcal{U}^{(1)}_{k,i} \cup\{(\beta_{k,i},u_{\beta_{k,i}},v_{\beta_{k,i}})\}|\beta\leq0\}$ has at most countable connected components. Furthermore, $S_0:=\{(\beta,u,v)\in\mathcal{U}^{(1)}_{k,i}|\beta=0\}$ has at most countable points. We divide the argument into two cases.

~~

\noindent{\bf Case 1. }There is a connected component $U_j$ of $\{(\beta,u,v)\in\mathcal{U}^{(1)}_{k,i} \cup\{(\beta_{k,i},u_{\beta_{k,i}},v_{\beta_{k,i}})\}|\beta\leq0\}$  such that $S_0\cap U_j\neq\emptyset$ and $(\beta_{k,i},u_{\beta_{k,i}},v_{\beta_{k,i}})\in U_j$. In this case, for any $(\beta,u,v)\in S_0\cap U_j$, $n(u)\neq k-1$ or $n(v)\neq k-1$, then $\mathcal{U}^{(1)}_{k,i}$ intersects with a semi-trivial branch for some $\beta\leq0$. This is impossible since the semi-trivial branch only bifurcate when $\beta>0$.

This implies that
$\mbox{Proj}_\beta\Big( \mathcal{U}^{(1)}_{k,i}\cup\{(\beta_{k,i},u_{\beta_{k,i}},v_{\beta_{k,i}})\} \Big)\subset(-\infty,0)$.

~~

\noindent{\bf Case 2. }There is a connected component $U_j$ of $\mathcal{U}^{(1)}_{k,i}\cup\{(\beta_{k,i},u_{\beta_{k,i}},v_{\beta_{k,i}})\}$  such that $S_0\cap U_j=\emptyset$ and $(\beta_{k,i},u_{\beta_{k,i}},v_{\beta_{k,i}})\in U_j$. Then, $\mbox{Proj}_\beta(U_j)\subset(-\infty,0)$. Then, we get
\begin{itemize}
  \item [$(1).$] $\mathcal{U}^{(1)}_{k,i} \cup\{(\beta_{k,i},u_{\beta_{k,i}},v_{\beta_{k,i}})\} =U_j\bigcup\Big( \mathcal{U}^{(1)}_{k,i}\cup\{(\beta_{k,i},u_{\beta_{k,i}},v_{\beta_{k,i}})\}\backslash U_j\Big)$;
  \item [$(2).$] both of $U_j$ and $\mathcal{U}^{(1)}_{k,i}\cup\{(\beta_{k,i},u_{\beta_{k,i}},v_{\beta_{k,i}})\}\backslash U_j$ are nonempty and closed;
  \item [$(3).$] $U_j\bigcap\Big( \mathcal{U}^{(1)}_{k,i}\cup\{(\beta_{k,i},u_{\beta_{k,i}},v_{\beta_{k,i}})\}\backslash U_j\Big)=\emptyset$.
\end{itemize}
This contradicts with the connectedness of $\mathcal{U}^{(1)}_{k,i}\cup\{(\beta_{k,i},u_{\beta_{k,i}},v_{\beta_{k,i}})\}$. Therefore, $\mathcal{U}^{(1)}_{k,i}$ does not intersect with any semi-trivial branch and Lemma \ref{l:Sign1} holds.
\begin{flushright}
$\Box$
\end{flushright}

\begin{lemma}\label{l:Boundedness}
Assume $i\neq k$.
For any compact interval $[\beta_*,\beta'_*]\subset\mathbb{R}$, there exists a constant $C>0$ such that of $(\beta,u,v)\in\mathcal{U}^{(1)}_{k,i}$ with $\beta\in[\beta_*,\beta'_*]$, we have $\|u\|,\|v\|\leq C$.
\end{lemma}

\noindent{\bf Proof.}
In this part, we will make use of Liouville-type Theorems \ref{t:Liouville} and \ref{t:Liouville2}. First, due to the above subsection, for any $(\beta,u,v)\in\mathcal{U}^{(1)}_{k,i}$, one has that $n(u)=n(v)=n(u+v)= k-1$ and $n(u-v)=i-1$. Notice we only need to prove the $L^\infty$-boundedness of $(u,v)\in\mathcal{U}^{(1)}_{k,i}$.

Arguing by contradiction, if there is a sequence $\beta_{n}\in [\beta_*,\beta'_*]$ such that there is a sequence of solutions $(u_n,v_n)$ to Problem (\ref{e:002}) with $\beta=\beta_{n}$ with
\begin{itemize}
  \item $n(u_n + v_n)= k-1$ and $n(u_n - v_n)= i-1$;
  \item for any $n$, there is a $x_n\in B_1$ such that $|v_n|_\infty\leq|u_n|_\infty=u_n(x_n)\to+\infty$;
  \item $\beta_{n}\to\beta_\infty\in[\beta_*,\beta'_*]$.
\end{itemize}

If $\beta_\infty\neq-1$,
define $M_n:=|u_n|_\infty$ and the re-scaled functions as
\begin{align}
\widetilde{u}_n=\frac{1}{M_n}u_n\Big(\frac{x+M_n x_n}{M_n}\Big)\nonumber
\end{align}
and
\begin{align}
\widetilde{v}_n=\frac{1}{M_n}v_n\Big(\frac{x+M_n x_n}{M_n}\Big).\nonumber
\end{align}
We get
\begin{equation}\label{e:scaled}
    \left\{
   \begin{array}{lr}
     -{\Delta}\widetilde{u}_n+\frac{1}{M_n^2}\widetilde{u}_n = (\widetilde{u}_n)^3+\beta_n (\widetilde{u}_n)(\widetilde{v}_n)^2\mbox{ in }\Omega_n ,\\
     -{\Delta}\widetilde{v}_n+\frac{1}{M_n^2}\widetilde{v}_n = (\widetilde{v}_n)^3+\beta_n (\widetilde{u}_n)^2(\widetilde{v}_n)\mbox{ in }\Omega_n ,\\
     \widetilde{u}_n,\widetilde{v}_n\in H_{0}^1(\Omega_n)
   \end{array}
   \right.
\end{equation}
with $\Omega_n:=M_n(B_1 -\{x_n\})$. Then, letting $n\to\infty$ and $\widetilde{u}_n\to U$ and $\widetilde{v}_n\to V$ in $C_{loc}^{2,\alpha}$. Due to a routine computation (cf. \cite{GidasSpruck1981,BartschDancerWang2010,Quittner2021,QuittnerBook,LiMiyagaki}), the solutions to Problem (\ref{e:scaled}) converge to the solution of the one of the following three problems.
\begin{equation}\label{e:scaled2}
    \left\{
   \begin{array}{lr}
     -{\Delta}U= U^3+\beta_\infty UV^2\mbox{ in }\mathbb{R}^3 ,\\
     -{\Delta}V= V^3+\beta_\infty U^2V\mbox{ in }\mathbb{R}^3 ,\\
     U,V\mbox{ are radial and in }C^{2,\alpha}(\mathbb{R}^3),\\
     n(U-V)\leq 2k+1\mbox{ and }n(U-V)\leq 2i+1;
   \end{array}
   \right.
\end{equation}
\begin{equation}\label{e:scaled3}
    \left\{
   \begin{array}{lr}
     -{\Delta}U= U^3+\beta_\infty UV^2\mbox{ in }\mathbb{R} ,\\
     -{\Delta}V= V^3+\beta_\infty U^2V\mbox{ in }\mathbb{R} ,\\
     U,V\mbox{ are even and in }C^{2,\alpha}(\mathbb{R}),\\
     n(U-V)\leq 2k+1\mbox{ and }n(U-V)\leq 2i+1;
   \end{array}
   \right.
\end{equation}
\begin{equation}\label{e:scaled4}
    \left\{
   \begin{array}{lr}
     -{\Delta}U= U^3+\beta_\infty UV^2\mbox{ in }\mathbb{R}_+ ,\\
     -{\Delta}V= V^3+\beta_\infty U^2V\mbox{ in }\mathbb{R}_+ ,\\
     U,V\mbox{ are in }C^{2,\alpha}(\mathbb{R}_+)\mbox{ and }U(0)=V(0)=0,\\
     n(U-V)\leq 2k+1\mbox{ and }n(U-V)\leq 2i+1.
   \end{array}
   \right.
\end{equation}
In
each of the cases, $U\neq0$ since $U(0)=\lim_{n\to\infty}\widetilde{u}_n(0)=1$. Then, we have a contradiction with Theorem \ref{t:Liouville} and Theorem \ref{t:Liouville2}.

If $\beta_\infty=-1$ and the limit equation is Problem (\ref{e:scaled2}), we have $U= V\equiv1$. In the following, we eliminate the case $U= V\equiv1$. If so, define $w_n=\widetilde{u}_n - \widetilde{v}_n$. Notice that for any bounded $\Omega\subset\Omega_n$, we have $\widetilde{u}_n,\widetilde{v}_n\to 1$ in $\Omega$ uniformly. Then, we have
\begin{equation}\label{e:scaled5}
    \left\{
   \begin{array}{lr}
     -{\Delta}w_n+\frac{1}{M^2_n}w_n=\big( \widetilde{u}_n^3 +\widetilde{v}_n^3+(1-\beta_n)\widetilde{u}_n\widetilde{v}_n\big)w_n\mbox{ in }\Omega ,\\
     |w_n|\leq2\mbox{ in }\Omega.
   \end{array}
   \right.
\end{equation}
Furthermore, we have $\beta_n\to-1$.
Letting $n\to\infty$, we get
\begin{align}
-\frac{d^2}{dr^2}w_n-\frac{2}{r}\frac{d}{dr}w_n=(4+o_n(1))w_n\mbox{ in }\Omega.
\end{align}
Here, $o_n(1)$ is the infinitesimal in $L^\infty(\Omega)$ norm.
Without loss of generality, let $\Omega$ be a sufficiently large ball.
Using Sturm's comparison theorem (cf. \cite[pp. 2]{SwansonBook1968}), there are $2i+1$ nodes of $w_n$ in $\Omega$. This contradicts with the assumption. Therefore, $U$ and $V$ cannot be $1$ and again we have a contradiction. The case of $U$ and $V$ solves Problem (\ref{e:scaled3}) and Problem (\ref{e:scaled4}) are the same.
\begin{flushright}
$\hfill \Box$
\end{flushright}

\noindent{\bf Proof of Assertions $(2-3)$ of Theorem \ref{t:bifurcation}.}
$\mathcal{U}^{(1)}_{k,i}$ cannot return to $\mathcal{T}^1_k$. Otherwise, it will contains a $(\beta,u,v)$ with $n(u-v)\neq i$. This is impossible due to Lemma \ref{l:Sign1}. Therefore, using the alternative result in Corollary \ref{coro:globalbranch}, it is an unbounded connected set. Again, applying Lemma \ref{l:Sign1}, $\mathcal{U}^{(1)}_{k,i}\subset\mathcal{S}^1_k$. This proves Assertion $(2)$ of Theorem \ref{t:bifurcation}.

The validity of Assertion $(3)$ of Theorem \ref{t:bifurcation} is evident due to Lemma \ref{l:Boundedness} and the fact that $0\notin\mbox{Proj}_\beta(\mathcal{U}^{(1)}_{k,i})$, i.e.,
$(-\infty,\beta_{k,i}) \subset\mbox{Proj}_\beta(\mathcal{U}^{(1)}_{k,i})\subset(-\infty,0)$.
\begin{flushright}
$\hfill \Box$
\end{flushright}

\subsection{Proofs of Assertions $(4-5)$ of Theorem \ref{t:bifurcation} and of Assertions $(2-5)$ of Theorem \ref{t:bifurcation2}}\label{subsection:(5)}

Assertion $(4)$ of Theorem \ref{t:bifurcation} can be proved via a similar approach as in the last subsection. We only sketch it here.
First of all, we have the following results.
\begin{lemma}\label{l:Sign2}
For any $i=1,\cdots,k-1$ and any $(\beta,u,v)\in\mathcal{U}^{(1)}_{k,i}$, we have that $n(u)=n(v)=n(u+v)=k-1$, $n(u-v)=i-1$ and $u(0),v(0)>0$.
\end{lemma}

\noindent{\bf Proof.}
Using a similar argument as in Lemma \ref{l:Sign1}, we only need to consider the case when $\mathcal{U}^{(1)}_{k,i}$ intersects with a semi-trivial branch. We only sketch the proof here. For any $k'\in\mathbb{N}$ and any $l=1,2,3,4$, we denote the bifurcation parameters of $\mathcal{ST}_{k'}^l$ by
\begin{align}
\widetilde{\beta}_{k',1}<\widetilde{\beta}_{k',2} <\cdots<\widetilde{\beta}_{k',i}<\cdots.\nonumber
\end{align}
Recall that $\widetilde{r}_{k',i}$ are the radii in Lemma \ref{l:1-parameterSemiTrivial} for $\mathcal{ST}_{k'}^{(1)}$.
As we pointed out in Remark \ref{r:SemiBranchParameters}, the bifurcation parameters and corresponding radii of $\mathcal{ST}_k^2$ are also $\widetilde{\beta}_{k,i}$ and $\widetilde{r}_{k,i}$.
Define the following sets.
\begin{align}
\mathcal{A}_1' & :=\Big(\mathcal{U}^{(1)}_{k',i}\cup\{(\beta_{k,i}, u_{\beta_{k,i}},v_{\beta_{k,i}})\}\Big)\cap \Big(\cup_{i\in\mathbb{N}_+,k'\neq k}\overline{\mathcal{B}_{\frac{\widetilde{r}_{k',i}}{4}}}(\widetilde{\beta}_{k',i}, w_k,0)\cup \cup_{i\in\mathbb{N}_+,k'\neq k}\overline{\mathcal{B}_{\frac{\widetilde{r}_{k',i}}{4}}}(\widetilde{\beta}_{k',i}, 0,w_k)\Big),\nonumber\\
\mathcal{B}'&:=\overline{\mathcal{U}^{(l)}_{k',i} \cup\{(\beta_{k,i},u_{\beta_{k,i}},v_{\beta_{k,i}})\} \backslash\mathcal{A}_1'},\nonumber\\
\mathcal{A}_2'&:=\Big(\mathcal{U}^{(1)}_{k',i}\cup\{(\beta_{k,i},u_{\beta_{k,i}},v_{\beta_{k,i}}) \}\Big)\cap \Big(\cup_{i\in\mathbb{N}_+,k'\neq k}\overline{\mathcal{B}_{\frac{\widetilde{r}_{k',i}}{2}}}(\widetilde{\beta}_{k',i}, w_k,0)\cup \cup_{i\in\mathbb{N}_+,k'\neq k}\overline{\mathcal{B}_{\frac{\widetilde{r}_{k',i}}{2}}}(\widetilde{\beta}_{k',i}, 0,w_k)\Big).\nonumber
\end{align}

As in the proof of Lemma \ref{l:Sign1}, we can find a connected component $\mathcal{C}'$ of $\mathcal{B}'$ such that $(\beta_{k,i},u_{\beta_{k,i}},v_{\beta_{k,i}})\in\mathcal{C}'$ and $\mathcal{C}'\cap\mathcal{A}'_2\neq\emptyset$. Then, $\mathcal{C}'$ intersects with a semi-trivial branch. This contradicts with the construction of $\mathcal{C}'$.
\begin{flushright}
$\Box$
\end{flushright}

\begin{lemma}\label{l:beta=1u=Cv}
Assume $\beta=1$. For any solution $(u,v)$ to Problem (\ref{e:001}), there exists a constant $C$ such that $u=Cv$. Therefore, $n(u)=n(v)$.
\end{lemma}

\noindent{\bf Proof.}
Since $(u,v)$ is a solution to Problem (\ref{e:001}), it holds that $n(u),n(v)<\infty$.
Notice that $u$ and $v$ are two eigenfunctions to the eigenvalue problem
\begin{equation}
    \left\{
   \begin{array}{lr}
     -{\Delta}\phi+\phi=\lambda (u^2+v^2)\phi\mbox{ in }B_1 ,\nonumber\\
     \phi\in H_{0,r}^1(B_1)\nonumber
   \end{array}
   \right.
\end{equation}
with $\lambda=1$. By \cite[Theorem XIII.7.53 and Corollary 7.56]{DunfordSchwartzBook1988}, the corresponding eigenspace is one-dimensional. Therefore, there exists a constant $C$ such that $u=Cv$.

\begin{flushright}
$\Box$
\end{flushright}

We have the following immediate consequence.

\begin{corollary}\label{coro:u=Cw}
Assume $\beta=1$. For any solution $(u,v)$ to Problem (\ref{e:001}) with $n(u)=n(v)=k-1$, there exists a constant $\theta\in[0,2\pi)$ such that $(u,v)=(\cos\theta\cdot\ w_k, \sin\theta\cdot w_k)$. Here, $w_k$ is the solution to Problem (\ref{e:002}).
\end{corollary}
\noindent{\bf Proof.}
By Lemma \ref{l:beta=1u=Cv}, $u$ and $v$ satisfy
\begin{align}
-\Delta u+u=(C^{-2}+1)u^3\nonumber
\end{align}
and
\begin{align}
-\Delta v+v=(C^2+1)v^3.\nonumber
\end{align}
Here, $C$ is the constant in Lemma \ref{l:beta=1u=Cv}.
By the uniqueness of the solution to Problem (\ref{e:002}) (see \cite[Theorem 1.3]{Tanaka2016}), we get $(u,v)=\Big(\frac{C}{\sqrt{1+C^2}}w,\frac{1}{\sqrt{1+C^2}}w\Big)$.
\begin{flushright}
$\Box$
\end{flushright}

\begin{remark}
The above results are proved in \cite[Theorem 4.2]{DaiTianZhang2019} and \cite[Lemma 2.2]{ZhouWang2020} of the case $k=1$.
\end{remark}

\noindent{\bf Proof of Assertions $(4-5)$ Theorem \ref{t:bifurcation}.}
The proof of $(4)$ of Theorem follows immediately if we notice that for any $(\beta,u,v)\in\mathcal{U}^{(1)}_{k,i}$ we have that $n(u)=n(v)=n(u+v)=k-1$, $n(u-v)=i-1$ and $u(0),v(0)>0$ and Theorems \ref{t:Liouville} and \ref{t:Liouville2}. To be precise, we have that $[\beta_{k,i},+\infty)\subset \mbox{Proj}_\beta (\mathcal{U}^{(1)}_{k,i})\subset(0,\infty)$.

Now we claim that $\mbox{Proj}_\beta(\mathcal{U}^{(1)}_{k,i})\subset(1,\infty)$. Otherwise, let us assume that there exists a sequence $(\beta_l,u_l,v_l)$ such that $\beta_l\downarrow1$ $n(u_l)=n(v_l)=n(u_l+v_l)=k-1$ and $n(u_l-v_l)<k-1$. By Lemma \ref{l:Boundedness} and Corollary \ref{coro:u=Cw}, there exists a $\theta_0\in[0,2\pi)$ such that
\begin{align}
(u_l,v_l)\to(\cos\theta\cdot w,\sin\theta\cdot w)\mbox{ as }l\to\infty.\nonumber
\end{align}
However, since either $\cos\theta\cdot w-\sin\theta\cdot w=0$ or $n(\cos\theta\cdot w-\sin\theta\cdot w)=k$,
by a similar argument as in the proof of Lemma \ref{l:Sign1}, we get that at lease one of the following assertions holds:
\begin{itemize}
  \item $u_l\to0$ in $H_{0,r}^1(B_1)$;
  \item $v_l\to0$ in $H_{0,r}^1(B_1)$.
\end{itemize}
Notice that the above two assertions can not hold at the same time because of the isolation of the trivial solution $(0,0)$. Without loss of the generality, assume that $\lim_{l\to\infty}v_l=0$. Then, $\lim_{l\to\infty}(u_l,v_l)\in\mathcal{ST}_k^{1}$. However, this contradicts with Lemma \ref{l:1-parameterSemiTrivial}. Therefore $1\notin\mbox{Proj}_\beta(\mathcal{U}^{(1)}_{k,i})$. We conclude $\mbox{Proj}_\beta(\mathcal{U}^{(1)}_{k,i})\subset(1,\infty)$ by its connectedness.

Assertion $(5)$ of Theorem \ref{t:bifurcation} holds evidently. Especially, the intersection $\overline{\mathcal{U}^{(1)}_{k,k}}\cap\mathcal{T}^1_k= \{(\beta_{k,k},u_{\beta_{k,k}},v_{\beta_{k,k}})\}$ is due to the bifurcation from a simple eigenvalue.

\begin{flushright}
$\Box$
\end{flushright}

\noindent{\bf Proof of Assertions $(2-5)$ Theorem \ref{t:bifurcation2}.}
Assertion $(2)$ of Theorem \ref{t:bifurcation2} follows a similar approach as Assertion $(2)$ of Theorem \ref{t:bifurcation}. We only prove Assertion $(3)$ of Theorem \ref{t:bifurcation2} since Assertion $(4)$ is similar. Especially, for any $(\beta,u,v)\in \mathcal{W}_{k,i}^{(1)}$, we get $u(0)>|v(0)|$. This is an immediate consequence of the fact that for any $(\beta,u,v)\in\mathcal{U}^{(1)}_{k,i}$ and the mapping (\ref{map:T}).
Moreover, by the definition of the mapping \ref{map:T}, it is known that for $i\geq k+1$, we have that
\begin{itemize}
  \item [$(a).$] $\mathcal{W}^{(1)}_{k,i}=T_{1}(\{(\beta,u,v) \in\mathcal{U}^{(1)}_{k,i}|\beta>-1\})$;
  \item [$(b).$] $0\notin\mbox{Proj}_\beta(\mathcal{U}^{(1)}_{k,i})$.
\end{itemize}
This implies that $3\notin\mbox{Proj}_\beta(\mathcal{W}^{(1)}_{k,i})$.
Assertion $(3)$ of Theorem \ref{t:bifurcation2} follows form the unboundedness of $\mathcal{W}^{(1)}_{k,i}$.
Assertion $(5)$ of Theorem \ref{t:bifurcation2} can be verified via a direct computation.
\begin{flushright}
$\Box$
\end{flushright}

\section{Proof of Theorem \ref{t:nonexistencebeta=3}}\label{section:nonexistence}
In this section, we prove Theorem \ref{t:nonexistencebeta=3}. To this end, we need the following Liouville-type theorem.
\begin{theorem}\label{t:Liouvillebeta>0}
Suppose that $\beta>0$, then problems
\begin{equation}
    \left\{
   \begin{array}{lr}
     -{\Delta}u=u^3+\beta uv^2 \,\,\,\,\,\,\,\,  \mbox{in}\ \mathbb{R}^N ,\\
     -{\Delta}v=v^3+\beta u^2v \,\,\,\,\,\,\,\,  \mbox{in}\ \mathbb{R}^N ,\\
     u,v\in C_b^2(\mathbb{R}^N)\mbox{ and }u,v\mbox{ are radial},
   \end{array}
   \right.
\end{equation}
\begin{equation}
    \left\{
   \begin{array}{lr}
     -u''=u^3+\beta uv^2 \,\,\,\,\,\,\,\,  \mbox{in}\ \mathbb{R} ,\\
     -v''=v^3+\beta u^2v \,\,\,\,\,\,\,\,  \mbox{in}\ \mathbb{R} ,\\
     u,v\in C_b^2(\mathbb{R})
   \end{array}
   \right.
\end{equation}
and
\begin{equation}
    \left\{
   \begin{array}{lr}
     -u''=u^3+\beta uv^2 \,\,\,\,\,\,\,\,  \mbox{in}\ [0,\infty) ,\\
     -v''=v^3+\beta u^2v \,\,\,\,\,\,\,\,  \mbox{in}\ [0,\infty) ,\\
     u,v\in C_b^2([0,\infty))\mbox{ with }u'(0)=v'(0)=0
   \end{array}
   \right.
\end{equation}
admits no solutions $(u,v)$ with $n(u),n(v)<\infty$.
\end{theorem}
This is a special case of \cite[Theorems 1.3 $\And$ 1.4]{LiMiyagaki}.

\Vs
\noindent{\bf Proof of Theorem \ref{t:nonexistencebeta=3}.}
By Theorem \ref{t:comparison} and (\ref{e:comparison}), it is evident that Theorem \ref{t:nonexistencebeta=3} holds for $\beta=3$. To be precise, we get
\begin{claim}\label{c:PneqQ}
For any $P,Q\in\mathbb{N}$ with $P\neq Q$, Problem (\ref{e:001}) admits no solution $(u,v)$ such that $n(u)=P$ and $n(v)=Q$ if $\beta=3$.
\end{claim}

Now we prove Theorem \ref{t:nonexistencebeta=3} by contradiction. Assume that there exists a sequence $\beta_l\to3$ such that
\begin{itemize}
  \item $\beta_l\neq 3$;
  \item for any $l=1,2,\cdots$, there exists a solution $(u_l,v_l)$ to Problem (\ref{e:001}) with $n(u_l)=P$, $n(v_l)=Q$ and $P\neq Q$.
\end{itemize}
By Theorem \ref{t:Liouvillebeta>0} and a similar computation in the proof of Lemma \ref{l:Boundedness}, we find a positive constant $C>0$ such that $\|u_l\|,\|v_l\|\leq C$. Thanks for the $(PS)$ condition, there exists $(u_\infty,v_\infty)\in H_{0,r}^1(B_1)\times H_{0,r}^1(B_1)$ with $(u_l,v_l)\to (u_\infty,v_\infty)$ in $H_{0,r}^1(B_1)$.

Then, $n(u_\infty)\neq n(v_\infty)$ due to $\beta_l\to3$ and Claim\ref{c:PneqQ}. Then, without loss of generality, one find a $r_0\in[0,1]$ such that $u_\infty(r_0)=\frac{d}{dr}u_\infty(r_0)=0$. This implies that $u_\infty\equiv0$.

If $v_\infty\neq 0$, then $(\beta_l,u_l,v_l)$ converges to a semi-trivial branch. This is impossible since $\beta=3$ is not a bifurcation parameter for semi-trivial branch, cf. Theorem \ref{t:bifurcation2}. If $v_\infty\equiv0$, this contradicts with the isolation of the trivial solution. Therefore, the assumption fails and Theorem \ref{t:nonexistencebeta=3} holds.

\begin{flushright}
$\Box$
\end{flushright}

\section{Proof of Theorem \ref{t:Liouville}}\label{section:MorseLiouville}


\noindent{\bf Proof of Theorem \ref{t:Liouville}.}
We divide our proof into two cases.

\vs
\noindent{\bf Case 1. $\beta\in(-1,3)$.}
Denote $u - v=w_1$ and $u+v=w_2$. Since $n(w_1),n(w_2)<\infty$, we get
\begin{align}
-\Delta w_1 &=-\Delta(u - v)= u^3 -  v^3+\beta uv(v - u)\nonumber\\
&= w_1( u^2 +  v^2+(1-\beta)uv)=\frac{1+\beta}{4}w_1^3+\frac{3-\beta}{4}w_1w_2^2.\nonumber
\end{align}
From a similar approach, we get
\begin{align}
-\Delta w_2=\frac{1+\beta}{4}w_2^3+\frac{3-\beta}{4}w_1^2w_2.\nonumber
\end{align}
Since $\beta\in(-1,3)$, we have $\frac{\beta+1}{4},\frac{3-\beta}{4}>0$. Then, \cite[Theorem 1.3]{LiMiyagaki} implies that $w_1\equiv w_2\equiv0$. Theorem \ref{t:Liouville} holds for $\beta\in(-1,3)$.

\vs
\noindent{\bf Case 2. $\beta<-1$.}
With a similar computation, we get
\begin{equation}\label{e:1111}
    \left\{
   \begin{array}{lr}
     -\Delta w_1=\frac{1+\beta}{4}w_1^3+\frac{3-\beta}{4}w_1w_2^2\mbox{ in }\mathbb{R}^N ,\\
     -\Delta w_2=\frac{1+\beta}{4}w_2^3+\frac{3-\beta}{4}w_1^2w_2\mbox{ in }\mathbb{R}^N ,\\
     w_1,w_2\in C_b^2(\mathbb{R}^N)\mbox{ and }w_1,w_2\mbox{ are radial}.
   \end{array}
   \right.
\end{equation}
Here, $\frac{1+\beta}{4}<0$ and $\frac{3-\beta}{4}>0$.
If $w_1\equiv0$, then $u\equiv v$ in $\mathbb{R}^N$ and $-\Delta u=(1+\beta)u^3$ in $\mathbb{R}^N$. Then, \cite[Lemma 2]{Brezis1984} and $\beta<-1$ imply that $u\equiv v\equiv0$. $w_2\equiv 0$ will draw a similar conclusion. Therefore,
without loss of generality, we can assume that both of $w_1$ and $w_2$ are nonzero. Since $n(u+v),n(u-v)<\infty$, we can find a large number $R_0>0$ such that for any $|x|>R_0$, we have $w_1(x), w_2(x)>0$. Let us consider the following argument on $B_{R_0}^c$.

Firstly, we claim that
\begin{claim}\label{c:1}
Assume $\beta<-1$. There is a constant $R_1\geq R_0$ such that $w_1(x)-w_2(x)$ has a constant sign for any $|x|>R_1$.
\end{claim}
We argue by contradiction. Suppose that there is a sequence of radii's $R_2<R_3<\cdots$ such that
\begin{align}
(-1)^j\big(w_1(x)-w_2(x)\big)>0\mbox{ for any }R_{j}<|x|<R_{j+1},
\end{align}
and $\lim_{j\to\infty}R_j=+\infty$. Let $R_{j_0}>R_0$. Then, on $R_{j_0}<|x|<R_{j_0 +1}$, without loss of generality, let $w_1>w_2(>0)$. By a direct computation, we get
\begin{align}
-\Delta(w_1 - w_2)=\frac{1+\beta}{4}(w_1^3 - w_2^3)+\frac{3-\beta}{4}w_1 w_2(w_2 - w_1)\leq 0.
\end{align}
The strong maximum principle (cf. \cite[Theorem 4/ pp. 333]{EvansBook}) gives that $w_1\equiv w_2$. This implies a contradiction with the uniqueness of ODE. A similar argument can prove $w_1<w_2$. Then, Claim \ref{c:1} is proved.

Based on this, we have that all of the zeroes of $u$ and $v$ are bounded in $B_{R_1 +1}(0)$.
By \cite[Theorem 1]{Quittner2021}, we complete the proof.

\vs
\noindent{\bf Step 3. $\beta\geq3$.}
We apply an idea in \cite{DancerWeiWeth2010}.
Since $n(w_1),n(w_2)<\infty$. Then, without loss of generality, there exists a $R>0$ such that for any $|x|>R$, we have $w_1(x),w_2(x)>0$. It is obvious that by Young's inequality
\begin{align}
-\Delta(w_1+w_2)\geq w_1^3+ w_2^3\geq c_0(w_1 + w_2)^3
\end{align}
for some constant $c_0>0$.
The following claim will conclude a contradiction.
\begin{claim}\label{c:Liouville}
For $c,R>0$, the following problem
\begin{equation}\label{e:111}
    \left\{
   \begin{array}{lr}
     -u''-\frac{N-1}{r}u'\geq cu^3\mbox{ in }(R,\infty),\\
     u>0\mbox{ in }(R,\infty),\\
     u\mbox{ is bounded on }(R,\infty)
   \end{array}
   \right.
\end{equation}
admits no nontrivial solution for $N=1,2,3$.
\end{claim}
\noindent{\bf Proof of Claim \ref{c:Liouville}.}
The case $N=2,3$ is handled by \cite[iii) of Theorem 3.3]{BidautVeronPohozaev2001}. Now let us consider the case $N=1$.
It is evident that $u''<0$ for any $r>R$. Then, $u'$ does not increase in $(R,\infty)$.
We argue by contradiction with the help of $\lim_{r\to\infty}u'(r)$. If there is a point $R_0>R$ such that $u'(r)<(>)0$, this will contradict with the boundedness of $u$. Then, $\lim_{r\to\infty}u'(r)=0$. Thereby, $u'(r)>0$ in $(R,\infty)$. The "=" at some point will lead us to a contradiction. Then, we know that $\lim_{r\to\infty}u(r)$ exist and $\lim_{r\to\infty}u(r)=:2\delta_0>0$.

Then, we get $-u''\geq cu^3$. And we have for some $r_0$,
\begin{align}
u'(r)=u'(r_0)+\int_{r_0}^r u''(s)ds\leq u'(r_0)-c\delta_0^3(r-r_0)<0
\end{align}
for large $r>0$. This is a contradiction.

\Vs

\noindent{\bf Case 4. $\beta=-1$.}
By a similar computation as in Case 1, we get
\begin{equation}\label{e:2222}
    \left\{
   \begin{array}{lr}
     -\Delta w_1=w_1w_2^2\mbox{ in }\mathbb{R}^N ,\\
     -\Delta w_2=w_1^2w_2\mbox{ in }\mathbb{R}^N ,\\
     w_1,w_2\in C_b^2(\mathbb{R}^N)\mbox{ and }w_1,w_2\mbox{ are radial}.
   \end{array}
   \right.
\end{equation}
Following the idea as in Claim \ref{c:1}, the following result is evident.
\begin{claim}
Assume $\beta=-1$. There is a constant $R_2$ such that $w_1(x)-w_2(x)$ has a constant sign for any $|x|>R_2$.
\end{claim}
Therefore, without loss of generality, there exists a constant $R_3$ such that for any $|x|\geq R_3$, $w_1(x)\geq w_2(x)>0$. Then, for any $x\in\mathbb{R}^N$ with $|x|> R_3$, we get
\begin{align}
-\Delta w_2\geq w_2^3.\nonumber
\end{align}
By Claim \ref{c:Liouville}, $w_2\equiv 0$ in $\mathbb{R}^N$. Then, $w_1$ is a constant and the proof is completed.

\begin{flushright}
$\Box$
\end{flushright}

%

\section{Miscellanies: Asymptotics of $\mathcal{W}^{(1)}_{k,i}$ for $i\geq k+1$ and of $\mathcal{U}^{(1)}_{k,i}$ for $i=1,\cdots,k-1$ as $\beta\to+\infty$}\label{section:miscellanies}

In this part, we consider the asymptotical behaviour of the branches $\mathcal{W}^{(1)}_{k,i}$ for $i\geq k+1$ and $\mathcal{U}^{(1)}_{k,i}$ for $i=1,\cdots,k-1$ as $\beta\to+\infty$. Our main results are the following.
\begin{theorem}\label{t:decay}
Let $i\geq k+1$. The following results hold.
\begin{itemize}
  \item [$(1).$] There are constants $C_{k,i},B>0$ such that for any $(\beta,u,v)\in\mathcal{W}^{(1)}_{k,i}$ with $\beta>B$, $|u|_\infty,|v|_\infty\leq\frac{C_{k,i}}{\sqrt{\beta}}$;
  \item [$(2).$] For the solution $(\beta,u,v)$ solving Problem (\ref{e:001}), as $\beta\to+\infty$, the functions $\sqrt{1+\beta}u\to u_\infty$ and $\sqrt{1+\beta}v\to v_\infty$ for some $(u_\infty,v_\infty)\in H_{0,r}^1(B_1)\times H_{0,r}^1(B_1)$ with
  \begin{equation}
  \left\{
  \begin{array}{lr}
     -{\Delta}u_\infty+u_\infty= u_\infty v_\infty^2\mbox{ in }B_1 ,\nonumber\\
     -{\Delta}v_\infty+v_\infty= u_\infty^2 v_\infty\mbox{ in }B_1 ,\nonumber\\
     u,v\in H_{0,r}^1(B_1)\nonumber
  \end{array}
  \right.
  \end{equation}
  and $n(u_\infty)=n(u_\infty+v_\infty)=n(u_\infty-v_\infty)=k-1$ and $n(v_\infty)=i-1$.
\end{itemize}

\end{theorem}

\begin{theorem}\label{t:decay2}
Let $i=1,\cdots,k-1$. The following results hold.
\begin{itemize}
  \item [$(1).$] There are constants $C'_{k,i},B>0$ such that for any $(\beta,u,v)\in\mathcal{U}^{(1)}_{k,i}$ with $\beta>B'$, $|u|_\infty,|v|_\infty\leq\frac{C'_{k,i}}{\sqrt{\beta}}$;
  \item [$(2).$] For the solution $(\beta,u,v)$ solving Problem (\ref{e:001}), as $\beta\to+\infty$, the functions $\sqrt{1+\beta}u\to u_\infty$ and $\sqrt{1+\beta}v\to v_\infty$ for some $(u_\infty,v_\infty)\in H_{0,r}^1(B_1)\times H_{0,r}^1(B_1)$ with
  \begin{equation}
  \left\{
  \begin{array}{lr}
     -{\Delta}u_\infty+u_\infty= u_\infty v_\infty^2\mbox{ in }B_1 ,\nonumber\\
     -{\Delta}v_\infty+v_\infty= u_\infty^2 v_\infty\mbox{ in }B_1 ,\nonumber\\
     u,v\in H_{0,r}^1(B_1)\nonumber
  \end{array}
  \right.
  \end{equation}
  and $n(u_\infty)=n(v_\infty)=n(u_\infty+v_\infty)=k-1$ and $n(u_\infty - v_\infty)=i-1$.
\end{itemize}
\end{theorem}
Since the proofs of Theorems \ref{t:decay} and \ref{t:decay2} are similar. We only prove Theorem \ref{t:decay}.

\noindent{\bf Proof of Theorem \ref{t:decay}.}
The main tool is the mapping $T$ defined in (\ref{map:T}). Recall that for the solution $(u,v)$ to the coupled system
\begin{equation}
    \left\{
   \begin{array}{lr}
     -{\Delta}u+u= u^3+\beta uv^2\mbox{ in }B_1 ,\nonumber\\
     -{\Delta}v+v= v^3+\beta u^2v\mbox{ in }B_1 ,\nonumber\\
     u,v\in H_{0,r}^1(B_1),\nonumber
   \end{array}
   \right.
\end{equation}
we have that the vector-valued function $T_1(u,v)=(u',v') =\Big(\frac{\sqrt{1+\beta}}{2}(u+v),\frac{\sqrt{1+\beta}}{2}(u-v)\Big)$ solves the system
\begin{equation}
    \left\{
   \begin{array}{lr}
     -{\Delta}u'+u'= (u')^3+\frac{3-\beta}{1+\beta} u'(v')^2\mbox{ in }B_1 ,\nonumber\\
     -{\Delta}v'+v'= (v')^3+\frac{3-\beta}{1+\beta} (u')^2v'\mbox{ in }B_1 ,\nonumber\\
     u',v'\in H_{0,r}^1(B_1).\nonumber
   \end{array}
   \right.
\end{equation}
We have that $(\beta,u,v)\in\mathcal{W}^{(1)}_{k,i}$ if and only if $\Big(\frac{3-\beta}{1+\beta},\frac{\sqrt{1+\beta}}{2}(u+v), \frac{\sqrt{1+\beta}}{2}(u-v)\Big)\in\mathcal{U}^{(1)}_{k,i}$ and that $\beta\uparrow+\infty$ if and only if $\frac{3-\beta}{1+\beta}\downarrow-1$. By a similar approach as in Lemma \ref{l:Boundedness}, there is a constant $C'_{k,i}>0$ such that as $\frac{3-\beta}{1+\beta}\downarrow-1$,
\begin{align}
\Bigg|\frac{\sqrt{1+\beta}}{2}(u+v)\Bigg|_\infty, \Bigg|\frac{\sqrt{1+\beta}}{2}(u-v)\Bigg|_\infty \leq C'_{k,i}.\nonumber
\end{align}
Assertion $(1)$ follows immediately. On the other hand, it is easy to check that for $(u,v)$ solving Problem (\ref{e:001}), one has
\begin{equation}
    \left\{
   \begin{array}{lr}
     -{\Delta}\Big(\sqrt{1+\beta}u\Big)+\Big(\sqrt{1+\beta}u\Big)= \frac{1}{1+\beta} \Big(\sqrt{1+\beta}u\Big)^3+ \frac{\beta}{1+\beta} \Big(\sqrt{1+\beta}u\Big)\Big(\sqrt{1+\beta}v\Big)^2\mbox{ in }B_1 ,\nonumber\\
     -{\Delta}\Big(\sqrt{1+\beta}v\Big)+\Big(\sqrt{1+\beta}v\Big)= \frac{1}{1+\beta} \Big(\sqrt{1+\beta}v\Big)^3+ \frac{\beta}{1+\beta} \Big(\sqrt{1+\beta}u\Big)^2\Big(\sqrt{1+\beta}v\Big)\mbox{ in }B_1 .\nonumber
   \end{array}
   \right.
\end{equation}
It is easy to verify that there exists $u_\infty,v_\infty\in H_{0,r}^1(B_1)$ such that $\Big(\sqrt{1+\beta}u,\sqrt{1+\beta}v\Big)\to(u_\infty,v_\infty)$ in $H_{0,r}^1(B_1)\times H_{0,r}^1(B_1)$.
This is due to the compactness of $\{(\beta,u,v)\in\mathcal{U}_{k,i}|\beta\in[-1,-1-\varepsilon]\}$ for small $\varepsilon>0$.

Next we prove that $u_\infty,v_\infty\neq 0$ by proving the limit of $(\frac{1}{2}\sqrt{1+\beta}(u+v), \frac{1}{2}\sqrt{1+\beta}(u-v))$ as $\frac{3-\beta}{1+\beta}\downarrow-1$ does not tend to $\{(u,v)\in H_{0,r}^1(B_1)\times H_{0,r}^1(B_1)|u=v\mbox{ or }u=-v\}$.

\begin{claim}\label{c:beta-1growth}
There exist positive constants $\varepsilon_0$ and $C$ such that if $\beta\in[-1,-1+\varepsilon_0]$, if $u=v$ or $u=-v$, we have either $\|u\|=\|v\|=0$ or $\|u\|,\|v\|\geq\frac{C}{\sqrt{1+\beta}}$.
\end{claim}
Indeed, multiplying $u$ on the both side of the first equation of Problem (\ref{e:001}) and integrating over $\Omega$, we get
\begin{align}
\int|\nabla u|^2+u^2=\int u^4+\beta\int u^2v^2=\int u^4+\beta\int u^4\leq C(1+\beta)\|u\|^4.\nonumber
\end{align}
Then, Claim \ref{c:beta-1growth} holds evidently. Since $\Big(\frac{3-\beta}{1+\beta},\frac{\sqrt{1+\beta}}{2}(u+v), \frac{\sqrt{1+\beta}}{2}(u-v)\Big)\in\mathcal{U}^{(1)}_{k,i}$, which does not contain $(0,0)$ and is bounded in $\mathbb{R}\times H_{0,r}^1(B_1)\times H_{0,r}^1(B_1)$, we prove that $(\frac{1}{2}\sqrt{1+\beta}(u+v), \frac{1}{2}\sqrt{1+\beta}(u-v))$ as $\frac{3-\beta}{1+\beta}\downarrow-1$ does not tend to $\{(u,v)\in H_{0,r}^1(B_1)\times H_{0,r}^1(B_1)|u=v\mbox{ or }u=-v\}$. Therefore, Theorem \ref{t:decay} follows.


\begin{flushright}
$\Box$
\end{flushright}


\begin{remark}
This is a generalization of \cite[Theorem 0.6]{SatoWang2013} to the case of radial nodal solutions.
\end{remark}

\appendix

\section{Appendix}\label{APPENDIX}
In this appendix, we prove Proposition \ref{p:MorseIndex}. Before we go to the proof, let us introduce the notion of a bump of a continuous radial function.
\begin{definition}\label{def:Bump}
For any radial domain $\Omega$ and a continuous radial function $u:\Omega\to\mathbb{R}$ with $n(u)=k$, due to Definition \ref{def:NodalNumber}, there exists positive numbers $x_0,x_1,\cdots,x_k\in(0,+\infty)$ satisfying
\begin{itemize}
  \item $\partial B_{x_i}(0)\subset\Omega$ for any $i=0,\cdots,k$;
  \item $u(x)|_{|x|=x_{i-1}}\cdot u(x)|_{|x|=x_{i}}<0$ for any $i=1,\cdots,k$.
\end{itemize}
Denote
\begin{align}
u_1(x)&:=u(x)\cdot\chi_{\{|x|\leq x_1\}}(x)\cdot\chi_{\{u(x)\cdot\mbox{sgn}(u(x_0))>0\}}(x);\nonumber\\
u_i(x)&:=u(x)\cdot\chi_{\{x_{i-2}<|x|<x_i\}}(x)\cdot\chi_{\{u(x)\cdot\mbox{sgn}(u(x_{i-1}))>0\}}(x)\mbox{ for }i=2,\cdots,k;\nonumber\\
u_{k+1}(x)&:=u(x)\cdot\chi_{\{|x|>x_{l-1}\}}(x)\cdot\chi_{\{u(x)\cdot\mbox{sgn}(u(x_k))>0\}}(x).\nonumber
\end{align}
We call $u_i$ the $i$-th bump of $u$. Here $\chi_S$ is the characteristic function of the set $S$.
\end{definition}

We recall first the minimax arguments to obtain the solution $w_k$ to Problem (\ref{e:002}).
\begin{itemize}
  \item [$(1).$] Let $0=r_0<r_1<\cdots<r_{k-1}<r_k=1$ be $k-1$ positive numbers;
  \item [$(2).$] Let $A_1=B_{r_1}(0)\subset B_1$ and $A_i=\{x\in B_1|r_{i-1}<|x|<r_i\}$ for $i=2,\cdots,k$;
  \item [$(3).$] Let $\mathcal{N}_i=\{u\in H_{0,r}^1(A_i)\backslash\{0\}|\int_{A_i}|\nabla u|^2 +|u|^2= \int_{A_i}|u|^4\}$ for $i=1,\cdots,k$.
\end{itemize}
Define
\begin{align}
c(A_i)=\inf_{u\in\mathcal{N}_i}\frac{1}{2}\int_{A_i}|\nabla u|^2+|u|^2-\frac{1}{4}\int_{A_i}|u|^4\nonumber
\end{align}
for $i=1,\cdots,k$ and
\begin{align}
c=\inf_{0<r_1<\cdots<r_{k-1}<1}\sum_{i=1}^k c(A_i).\nonumber
\end{align}
\cite{BartschWillem1993} and \cite[Chapter 5]{WillemBook1996} ensure the validity of the following results.
\begin{lemma}
Under the above assumptions, we have
\begin{itemize}
  \item [$(1).$] There exist $k-1$ numbers $r^0_1<\cdots<r^0_{k-1}$ achieving $c$;
  \item [$(2).$] Denote $A^0_1=B_{r^0_1}(0)$ and $A^0_i=\{x\in\mathbb{R}^3|r^0_{i-1}<|x|<r^0_{i}\}$ for $i=2,\cdots,k$. For any $i=1,\cdots,k$, there is a positive function $u^0_i\in H_{0,r}^1(A^0_i)$ such that $c(A^0_i)=\frac{1}{2}\int_{A^0_i}|\nabla u^0_i|^2+|u^0_i|^2-\frac{1}{4}\int_{A^0_i}|u^0_i|^4$;
  \item [$(3).$] $u_k:=\sum_{i=1}^k(-1)^{i-1}u^0_i$ is a solution to Problem (\ref{e:002}) changing its sign exactly $k-1$ times.
\end{itemize}
\end{lemma}

Denote the Nehari set
\begin{align}
\mathcal{N}_k^*=\Bigg\{\sum_{i=1}^{k}(-1)^{i-1}u_i\Bigg| & \exists0=r_0<r_1<\cdots<r_{k-1}<r_{k}=1\mbox{ such that }\nonumber\\
&0\leq u_i\in H_{0,r}^1(A_i)\mbox{ for }i=1,\cdots,k\Bigg\}.\nonumber
\end{align}
It is evident that $c=\inf_{u\in\mathcal{N}_k^*}I(u)=I(w_k)$.
Now we begin to prove Proposition \ref{p:MorseIndex}.

\vs

\noindent{\bf Proof of Proposition \ref{p:MorseIndex}.}
First, it is evident that $m(w_k)\geq k$. This is due to the fact that
\begin{align}
D^2 I(w_k)[w_{k,i},w_{k,i}]=\int|\nabla w_{k,i}|^2+|w_{k,i}|^2-3\int|w_{k,i}|^4=-2\int|w_{k,i}|^4<0\nonumber
\end{align}
for any $i=1,\cdots,k$. Here, $w_{k,i}$ denotes the $i$-th bump of $w_k$. See Definition \ref{def:Bump}.

It is sufficient to prove that $m(w_k)\leq k$. We argue by contradiction. Suppose that $m(u)\geq k+1$.
Denote
\begin{align}
D=\left\{\sum_{i=1}^k t_i w_{k,i}\Bigg|t_i>0\mbox{ for }i=1,\cdots,k\right\}.\nonumber
\end{align}
Define a mapping as follows
\begin{align}
L:(0,+\infty)^k&\to D\nonumber\\
(t_1,\cdots,t_k)&\mapsto\sum_{i=1}^k t_i w_{k,i}.\nonumber
\end{align}
It is easy to verify that $L$ is a homomorphism and $L(1,\cdots,1)=w_k$.
Denote $V^-$ the negative space of $D^2 I(w_k)$. On the other hand, since we assume that $m(u)\geq k+1$, there is a function $\phi_*\in V^-$ such that $\phi_*\notin\mbox{span}\{w_{k,1},\cdots,w_{k,k}\}$. Without loss of generality, we can assume that $\phi_*\in C^1(\overline{B_1})$. Define a cut-off function $\varphi\in C^\infty([0,+\infty),[0,1])$ such that
\begin{equation}
\varphi(s)=\left\{
\begin{aligned}
1 & \qquad & s\in[0,1], \nonumber\\
0 & \qquad & s\geq2.\nonumber
\end{aligned}
\right.
\end{equation}
For small constants $\delta_1,\delta_2>0$,
define a mapping $S: D\to H_{0,r}^1(B_1)$ as
\begin{align}
S(u)=u+\delta_1\cdot\phi_* \cdot \varphi\Big(\frac{|L^{-1}(u)-(1,\cdots,1)|}{\delta_2}\Big).\nonumber
\end{align}
Since $\phi_*\in C^1(\overline{B_1})$, there is a $\delta_1$ such that for any $u\in S(D)$, $u(0)>0$ and $n(u)=k-1$.

On one hand, using Morse lemma at $w_k$, cf. \cite[Lemma 4.1/pp. 33]{ChangBook1993}, there exist positive constants $\delta_3,\delta_4$ such that
\begin{align}\label{inclu:contradiction1}
S(D)\cap B_{\delta_3}(w_k)\subset I^{c-\delta_4}.
\end{align}
Here, $B_{\delta_3}(w_k)=\{u\in H_{0,r}^1(B_1)|\|u-w_k\|<\delta_3\}$ and $c=\inf_{u\in\mathcal{N}_k^*}I(u)=I(w_k)$. Notice that we can chose $\delta_1<\delta_3$.

On the other hand, consider the functional $\xi(u):=\int|\nabla u|^2+|u|^2-\int|u|^4$. Since for any $u\in S(D)$, $u(0)>0$ and $n(u)=k-1$, the following mapping is well-defined
\begin{align}
X(t_1,\cdots,t_k)=\Bigg(\xi\Big(\big(S\circ L(t_1,\cdots,t_k)\big)_1\Big),\cdots,\xi\Big(\big(S\circ L(t_1,\cdots,t_k)\big)_k\Big)\Bigg).\nonumber
\end{align}
Here, $\big(S\circ L(t_1,\cdots,t_k)\big)_i$ denotes the $i$-th bump of the radial function $S\circ L(t_1,\cdots,t_k)$.
Notice that $X:(0,+\infty)^k\to (0,+\infty)^k$ is a continuous mapping. In order to contradict with (\ref{inclu:contradiction1}), we need to find a $(t^0_1,\cdots,t_k^0)\in B_{\delta_1}((1,\cdots,1))$ with $S\circ L(t^0_1,\cdots,t^0_k)\in\mathcal{N}_k^*$. If so, $I\big(S\circ L(t^0_1,\cdots,t^0_k)\big)\geq c=\inf_{u\in\mathcal{N}_k^*}I(u)$. To this end, we need to check the following claim.
\begin{claim}\label{c:intersection}
There is a $(t^0_1,\cdots,t_k^0)\in B_{\delta_2}((1,\cdots,1))$ such that $X(t_1^0,\cdots,t_k^0)=(0,\cdots,0)$.
\end{claim}

\noindent{\bf Proof of Claim \ref{c:intersection}.}
Denote $X(t_1,\cdots,t_k)=(X_1(t_1,\cdots,t_k),\cdots, X_k(t_1,\cdots,t_k))$. Without loss of generality, we assume that the constants $\delta_2,\delta_3<\frac{1}{16}$. Then, for $i=1,\cdots,k$,
\begin{itemize}
  \item $X_i(t_1,\cdots,t_k)|_{t_i=\frac{1}{2}}=\frac{1}{4}\int|\nabla w_{k,i}|^2+|w_{k,i}|^2-\frac{1}{16}\int|w_{k,i}|^4>0$;
  \item $X_i(t_1,\cdots,t_k)|_{t_i=2}=4\int|\nabla w_{k,i}|^2+|w_{k,i}|^2-16\int|w_{k,i}|^4<0$.
\end{itemize}
The existence of $(t_1^0,\cdots,t_k^0)\in [\frac{1}{2},2]^k\subset[0,+\infty)^k$ such that $X(t_1^0,\cdots,t_k^0)=(0,\cdots,0)$ follows Poincar\'e-Miranda theorem, cf. \cite[pp. 547]{Kulpa1997}. $(t_1^0,\cdots,t_k^0)\in B_{\delta_2}((1,\cdots,1))$ is due to the fact that $X_i|_{B_{\delta_2}((1,\cdots,1))^c}=t_i-1$ for $i=1,\cdots,k$.
\begin{flushright}
$\Box$
\end{flushright}


\noindent{\bf Acknowledgment.}
The authors would like to express their sincere appreciation to Xiaopeng Huang for his technical assistance. HL was supported by FAPESP Proc 2022/15812-0. OHM was supported by CNPq Proc 303256/2022-2 and FAPESP Proc 2022/16407-1.




\vspace{-0.11cm}

{

}

\end{sloppypar}
\end{document}